\documentclass[11pt,reqno]{article}
\usepackage{amsfonts}
\usepackage{graphicx}

\begin{document}

\def\R{{\mathbb R}}
\def\T{{\mathbb T}}
\def\S{{\mathbb S}}
\def\C{{\mathbb C}}
\def\Z{{\mathbb Z}}
\def\N{{\mathbb N}}
\def\H{{\mathbb H}}
\def\B{{\mathbb B}}
\def\diam{\mbox{\rm diam}}
\def\sn{\S^{n-1}}
\def\rr{{\cal R}}
\def\mt{{\Lambda}}
\def\e{\emptyset}
\def\dQ{\partial Q}
\def\dk{\partial K}
\def\endofproof{{\rule{6pt}{6pt}}}
\def\di{\displaystyle}
\def\dist{\mbox{\rm dist}}
\def\sa+{\Sigma_A^+}
\def\du{\frac{\partial}{\partial u}}
\def\dv{\frac{\partial}{\partial v}}
\def\dt{\frac{d}{d t}}
\def\dx{\frac{\partial}{\partial x}}
\def\con{\mbox{\rm const }}
\def\nn{{\cal N}}
\def\mm{{\cal M}}
\def\kk{{\cal K}}
\def\ll{{\cal L}}
\def\vv{{\cal V}}
\def\bb{{\cal B}}
\def\ma{\mm_{a}}
\def\lab{L_{ab}}
\def\mabn{\mm_{a}^N}
\def\man{\mm_a^N}
\def\labn{L_{ab}^N}
\def\fa{f^{(a)}}
\def\ff{{\cal F}}
\def\i{{\bf i}}
\def\gge{{\cal G}_\epsilon}
\def\gej{\chi^{(j)}_\mu}
\def\ge{\chi_\epsilon}
\def\chio{\chi^{(1)}}
\def\chit{\chi^{(2)}}
\def\chii{\chi^{(i)}}
\def\chil{\chi^{(\ell)}}
\def\gett{\chi^{(2)}_{\mu}}
\def\geol{\chi^{(1)}_{\ell}}
\def\getl{\chi^{(2)}_{\ell}}
\def\geil{\chi^{(i)}_{\ell}}
\def\gee{\chi_{\ell}}
\def\tt{{\cal T}}
\def\uu{{\cal U}}
\def\wloc{W_{\epsilon}}
\def\Int{\mbox{\rm Int}}
\def\dist{\mbox{\rm dist}}
\def\pr{\mbox{\rm pr}}
\def\pp{{\cal P}}
\def\tpp{\widetilde{\pp}}
\def\aa{{\cal A}}
\def\cc{{\cal C}}
\def\supp{\mbox{\rm supp}}
\def\Arg{\mbox{\rm Arg}}
\def\In{\mbox{\rm Int}}
\def\con{\mbox{\rm const}\;}
\def\Re{\mbox{\rm Re}}
\def\li{\mbox{\rm li}} 
\def\Seo{S^*_\epsilon(\Omega)}
\def\sdk{S^*_{\dk}(\Omega)}
\def\lae{\Lambda_{\epsilon}}
\def\ep{\epsilon}
\def\oo{{\cal O}}
\def\be{\begin{equation}}
\def\ee{\end{equation}}
\def\beqn{\begin{eqnarray*}}
\def\eeqn{\end{eqnarray*}}
\def\Pr{\mbox{\rm Pr}}

\def\gi{\gamma^{(i)}}
\def\ii{{\imath }}
\def\jj{{\jmath }}
\def\II{{\cal I}}
\def\ccij{ \cc_{i'_0,j'_0}[\eta]}
\def\dd{{\cal D}}
\def\la{\langle}
\def\ra{\rangle}
\def\bs{\bigskip}
\def\xio{\xi^{(0)}}
\def\xo{x^{(0)}}
\def\zo{z^{(0)}}
\def\Con{\mbox{\rm Const}\;}
\def\do{\partial \Omega}
\def\dk{\partial K}
\def\dl{\partial L}
\def\ll{{\cal L}}
\def\kk{{\cal K}}
\def\kk{{\cal K}}
\def\pr{{\rm pr}}
\def\ff{{\cal F}}
\def\G{{\cal G}}
\def\C{{\bf C}}
\def\dist{{\rm dist}}
\def\dds{\frac{d}{ds}}
\def\con{{\rm const}\;}
\def\Con{{\rm Const}\;}
\def\di{\displaystyle}
\def\oo{\mbox{\rm O}}
\def\hess{\mbox{\rm Hess}}
\def\gi{\gamma^{(i)}}
\def\endofproof{{\rule{6pt}{6pt}}}
\def\xm{x^{(m)}}
\def\vm{\varphi^{(m)}}
\def\km{k^{(m)}}
\def\dm{d^{(m)}}
\def\kam{\kappa^{(m)}}
\def\dem{\delta^{(m)}}
\def\xim{\xi^{(m)}}
\def\ep{\epsilon}
\def\ms{\medskip}
\def\ex{\mbox{\rm extd}}

\def\clip{C^{\mbox{\footnotesize \rm Lip}}}
\def\wlocs{W^s_{\mbox{\footnote\rm loc}}}
\def\Lip{\mbox{\rm Lip}}

\def\Xr{X^{(r)}}
\def\lip{\mbox{{\footnotesize\rm Lip}}}
\def\Vol{\mbox{\rm Vol}}

\def\naf{\nabla f(z)}
\def\so{\sigma_0}
\def\Xo{X^{(0)}}
\def\z1{z^{(1)}}
\def\Vo{V^{(0)}}
\def\Yo{Y{(0)}}

\def\uo{u^{(0)}}
\def\vo{v^{(0)}}
\def\no{\nu^{(0)}}
\def\psa{\partial^{(s)}_a}
\def\hcd{\hc^{(\delta)}}
\def\Md{M^{(\delta)}}
\def\Uo{U^{(1)}}
\def\Ut{U^{(2)}}
\def\Uj{U^{(j)}}
\def\no{n^{(1)}}
\def\nt{n^{(2)}}
\def\nj{n^{(j)}}
\def\ccm{\cc^{(m)}}

\def\ooo{\oo^{(1)}}
\def\oot{\oo^{(2)}}
\def\ooj{\oo^{(j)}}
\def\fo{f^{(1)}}
\def\ft{f^{(2)}}
\def\fj{f^{(j)}}
\def\wo{w^{(1)}}
\def\wt{w^{(2)}}
\def\wj{w^{(j)}}
\def\Vo{V^{(1)}}
\def\Vt{V^{(2)}}
\def\Vj{V^{(j)}}

\def\Ul{U^{(\ell)}}
\def\Uj{U^{(j)}}
\def\wl{w^{(\ell)}}
\def\Vl{V^{(\ell)}}
\def\Ujj{U^{(j+1)}}
\def\wjj{w^{(j+1)}}
\def\Vjj{V^{(j+1)}}
\def\Ujo{U^{(j_0)}}
\def\wjo{w^{(j_0)}}
\def\Vjo{V^{(j_0)}}
\def\vj{v^{(j)}}
\def\vl{v^{(\ell)}}

\def\f0{f^{(0)}}

\def\gl{\gamma_\ell}
\def\id{\mbox{\rm id}}
\def\piU{\pi^{(U)}}

\def\cca{C^{(a)}}
\def\bba{B^{(a)}}
\def\co{\; \stackrel{\circ}{C}}
\def\Lo{\; \stackrel{\circ}{L}}

\def\oV{\overline{V}}
\def\saa{\Sigma^+_A}
\def\sa{\Sigma_A}
\def\mta{\Lambda(A, \tau)}
\def\mtaa{\Lambda^+(A, \tau)}

\def\Int{\mbox{\rm Int}}
\def\epo{\ep^{(0)}}
\def\pH{\partial \H^{n+1}}
\def\sh{S^*(\H^{n+1})}
\def\zoo{z^{(1)}}
\def\yoo{y^{(1)}}
\def\xoo{x^{(1)}}


\def\supp{\mbox{\rm supp}}
\def\Arg{\mbox{\rm Arg}}
\def\In{\mbox{\rm Int}}
\def\diam{\mbox{\rm diam}}
\def\e{\emptyset}
\def\endofproof{{\rule{6pt}{6pt}}}
\def\di{\displaystyle}
\def\dist{\mbox{\rm dist}}
\def\con{\mbox{\rm const }}
\def\Box{\spadesuit}
\def\Int{\mbox{\rm Int}}
\def\dist{\mbox{\rm dist}}
\def\pr{\mbox{\rm pr}}
\def\be{\begin{equation}}
\def\ee{\end{equation}}
\def\beqn{\begin{eqnarray*}}
\def\eeqn{\end{eqnarray*}}
\def\la{\langle}
\def\ra{\rangle}
\def\bs{\bigskip}
\def\Con{\mbox{\rm Const}\;}
\def\clip{C^{\alpha}}
\def\wlocs{W^s_{\mbox{\footnote\rm loc}}}
\def\Lip{\mbox{\rm Lip}}
\def\lip{\mbox{\footnotesize\rm Lip}}
\def\Re{\mbox{\rm Re}}
\def\li{\mbox{\rm li}} 
\def\ep{\epsilon}
\def\ms{\medskip}
\def\dds{\frac{d}{ds}}
\def\oo{\mbox{\rm O}}
\def\hess{\mbox{\rm Hess}}
\def\id{\mbox{\rm id}}
\def\ii{{\imath }}
\def\jj{{\jmath }}
\def\graph{\mbox{\rm graph}}
\def\span{\mbox{\rm span}}
\def\Intu{\Int^u}
\def\Ints{\Int^s}

\def\i{{\bf i}}
\def\C{{\bf C}}

\def\ss{{\cal S}}
\def\tt{{\cal T}}
\def\E{{\cal E}}
\def\rr{{\cal R}}
\def\nn{{\cal N}}
\def\mm{{\cal M}}
\def\kk{{\cal K}}
\def\ll{{\cal L}}
\def\vv{{\cal V}}
\def\ff{{\cal F}}
\def\hh{{\cal H}}
\def\tt{{\cal T}}
\def\uu{{\cal U}}
\def\cc{{\cal C}}
\def\pp{{\cal P}}
\def\aa{{\cal A}}
\def\oo{{\cal O}}
\def\II{{\cal I}}
\def\dd{{\cal D}}
\def\ll{{\cal L}}
\def\ff{{\cal F}}
\def\G{{\cal G}}

\def\hs{\hat{s}}
\def\hz{\hat{z}}
\def\hL{\widehat{L}}
\def\hl{\hat{l}}
\def\hl{\hat{l}}
\def\hc{\hat{\cc}}
\def\hbb{\widehat{\cal B}}
\def\hu{\hat{u}}
\def\hX{\hat{X}}
\def\hx{\hat{x}}
\def\hu{\hat{u}}
\def\hv{\hat{v}}
\def\hQ{\hat{Q}}
\def\hC{\widehat{C}}
\def\hF{\hat{F}}
\def\hf{\hat{f}}
\def\hii{\hat{\ii}}
\def\hr{\hat{r}}
\def\hq{\hat{q}}
\def\hy{\hat{y}}
\def\hZ{\widehat{Z}}
\def\hz{\hat{z}}
\def\hE{\widehat{E}}
\def\hR{\widehat{R}}
\def\hell{\hat{\ell}}
\def\hs{\hat{s}}
\def\hW{\widehat{W}}
\def\hS{\widehat{S}}
\def\hV{\widehat{V}}
\def\hB{\widehat{B}}
\def\hhh{\widehat{\cal H}}
\def\hK{\widehat{K}}
\def\hU{\widehat{U}}
\def\hhh{\widehat{\hh}}
\def\hdd{\widehat{\dd}}
\def\hZ{\widehat{Z}}
\def\hGa{\widehat{\Gamma}}

\def\hal{\hat{\alpha}}
\def\hbe{\hat{\beta}}
\def\hg{\hat{\gamma}}
\def\hrho{\hat{\rho}}
\def\hd{\hat{\delta}}
\def\hphi{\hat{\phi}}
\def\hmu{\hat{\mu}}
\def\hnu{\hat{\nu}}
\def\hsi{\hat{\sigma}}
\def\htau{\hat{\tau}}
\def\hpi{\hat{\pi}}
\def\hep{\hat{\epsilon}}
\def\hxi{\hat{\xi}}
\def\hLa{\widehat{\Lambda}^u}
\def\hPhi{\widehat{\Phi}}
\def\hPsi{\widehat{\Psi}}
\def\hPhii{\widehat{\Phi}^{(i)}}
\def\hath{\hat{h}}

\def\tc{\tilde{C}}
\def\tg{\tilde{\gamma}}  
\def\tV{\widetilde{V}}
\def\tC{\widetilde{\cc}}
\def\tR{\widetilde{R}}
\def\tb{\tilde{b}}
\def\tt{\tilde{t}}
\def\tx{\tilde{x}}
\def\tp{\tilde{p}}
\def\tz{\tilde{Z}}
\def\tZ{\tilde{Z}}
\def\tF{\tilde{F}}
\def\tf{\tilde{f}}
\def\tp{\tilde{p}}
\def\te{\tilde{e}}
\def\tv{\tilde{v}}
\def\tu{\tilde{u}}
\def\tw{\tilde{w}}
\def\ts{\tilde{\sigma}}
\def\tr{\tilde{r}}
\def\tU{\widetilde{U}}
\def\tS{\tilde{S}}
\def\tP{\widetilde{P}}
\def\ttau{\tilde{\tau}}
\def\tLip{\widetilde{\Lip}}
\def\tz{\tilde{z}}
\def\tS{\tilde{S}}
\def\tts{\tilde{\sigma}}
\def\tVl{\widetilde{V}^{(\ell)}}
\def\tVj{\widetilde{V}^{(j)}}
\def\tVo{\widetilde{V}^{(1)}}
\def\tVj{\widetilde{V}^{(j)}}
\def\tPsi{\tilde{\Psi}}
 \def\tp{\tilde{p}}
 \def\tVjo{\widetilde{V}^{(j_0)}}
\def\tvj{\tilde{v}^{(j)}}
\def\tVjj{\widetilde{V}^{(j+1)}}
\def\tvl{\tilde{v}^{(\ell)}}
\def\tVt{\widetilde{V}^{(2)}}
\def\tR{\widetilde{R}}
\def\tQ{\widetilde{Q}}
\def\oL{\tilde{\Lambda}}
\def\tq{\tilde{q}}
\def\tx{\tilde{x}}
\def\ty{\tilde{y}}
\def\tz{\tilde{z}}
\def\txo{\tilde{x}^{(0)}}
\def\tso{\tilde{\sigma}_0}
\def\tmt{\tilde{\Lambda}}
\def\tg{\tilde{g}}
\def\tsi{\tilde{\sigma}}
\def\tC{\tilde{C}}
\def\tc{\tilde{c}}
\def\tell{\tilde{\ell}}
\def\trho{\tilde{\rho}}
\def\ts{\tilde{s}}
\def\tB{\widetilde{B}}
\def\thh{\widetilde{\cal H}}
\def\tV{\widetilde{V}}
\def\trr{\tilde{r}}
\def\te{\tilde{e}}
\def\tv{\tilde{v}}
\def\tu{\tilde{u}}
\def\tw{\tilde{w}}
\def\trho{\tilde{\rho}}
\def\tell{\tilde{\ell}}
\def\tz{\tilde{Z}}
\def\tF{\tilde{F}}
\def\tf{\tilde{f}}
\def\tp{\tilde{p}}
\def\ttau{\tilde{\tau}}
\def\tz{\tilde{z}}
\def\tg{\tilde{\gamma}}  
\def\tV{\widetilde{V}}
\def\tC{\widetilde{\cc}}
\def\tLa{\widetilde{\Lambda}^u}
\def\tR{\widetilde{R}}
\def\tr{\tilde{r}}
\def\tc{\widetilde{C}}
\def\tD{\widetilde{D}}
\def\tt{\tilde{t}}
\def\tp{\tilde{p}}
\def\tS{\tilde{S}}
\def\tts{\tilde{\sigma}}
\def\tZ{\widetilde{Z}}
\def\tdelta{\tilde{\delta}}
\def\th{\tilde{h}}
\def\tB{\widetilde{B}}
\def\thh{\widetilde{\hh}}
\def\tep{\tilde{\ep}}
\def\tE{\widetilde{E}}
\def\tu{\tilde{u}}
\def\txi{\tilde{\xi}}
\def\teta{\tilde{\eta}}
\def\tRR{\widetilde{\rr}}

\def\sr{{\sc r}}
\def\mt{{\Lambda}}
\def\do{\partial \Omega}
\def\dk{\partial K}
\def\dl{\partial L}
\def\wloc{W_{\epsilon}}
\def\piU{\pi^{(U)}}
\def\Rio{\R_{i_0}}
\def\Ri{\R_{i}}
\def\Rii{\R^{(i)}}
\def\Riii{\R^{(i-1)}}
\def\hRii{\widehat{\R}_i}
\def\hRiio{\widehat{\R}_{(i_0)}}
\def\Eii{E^{(i)}}
\def\Eio{E^{(i_0)}}
\def\Rj{\R_{j}}
\def\Vio{{\cal V}^{i_0}}
\def\Vi{{\cal V}^{i}}
\def\Wio{W^{i_0}}
\def\Wioo{W^{i_0-1}}
\def\hi{h^{(i)}}
\def\Psii{\Psi^{(i)}}
\def\pii{\pi^{(i)}}
\def\piii{\pi^{(i-1)}}
\def\gxyii{g_{x,y}^{i-1}}
\def\span{\mbox{\rm span}}
\def\Jac{\mbox{\rm Jac}}
\def\Vol{\mbox{\rm Vol}}
\def\limp{\lim_{p\to\infty}}
\def\hh{{\mathcal H}}

\def\yijl{Y_{i,j}^{(\ell)}}
\def\xijl{X_{i,j}^{(\ell)}}
\def\xij{X_{i,j}}
\def\hyijl{\widehat{Y}_{i,j}^{(\ell)}}
\def\hxijl{\widehat{X}_{i,j}^{(\ell)}}
\def\hxij{\widehat{X}_{i,j}}
\def\eijl{\omega_{i,j}^{(\ell)}}
\def\eij{\omega_{i,j}}
\def\Gl{\Gamma_\ell}

\def\hLao{\widehat{\Lambda}^{u,1}}
\def\tLao{\widetilde{\Lambda}^{u,1}}
\def\Lao{\Lambda^{u,1}}
\def\cLao{\check{\Lambda}^{u,1}}
\def\cB{\check{B}}
\def\tpi{\tilde{\pi}}
\def\J{{\sf J}}
\def\bJ{{\mathbb J}}

\def\hcc{\widehat{\cc}}
\def\hpp{\widehat{\pp}}
\def\ttP{\widetilde{\pp}}
\def\tP{\widetilde{P}}
\def\hP{\widehat{P}}
\def\hY{\widehat{Y}}

\def\diamtef{{\footnotesize \diam_\theta}}


\def\tc{\tilde{C}}
\def\tg{\tilde{\gamma}}  
\def\tV{\widetilde{V}}
\def\tW{\widetilde{W}}
\def\tC{\widetilde{\cc}}
\def\tKo{\widetilde{K_0}}
\def\tUKo{\widetilde{U\setminus K_0}}

\def\wo{w^{(1)}}
\def\vo{v^{(1)}}
\def\uo{u^{(1)}}
\def\wt{w^{(2)}}
\def\xio{\xi^{(1)}}
\def\xit{\xi^{(2)}}
\def\etao{\eta^{(1)}}
\def\etat{\eta^{(2)}}
\def\zetao{\zeta^{(1)}}
\def\zetat{\zeta^{(2)}}
\def\vt{v^{(2)}}
\def\ut{u^{(2)}}
\def\Wo{W^{(1)}}
\def\Vo{V^{(1)}}
\def\Uo{U^{(1)}}
\def\Wt{W^{(2)}}
\def\Vt{V^{(2)}}
\def\Ut{U^{(2)}}
\def\tmu{\tilde{\mu}}
\def\tla{\tilde{\lambda}}
\def\diamf{{\rm\footnotesize diam}}
\def\Intu{\mbox{\rm Int}^u}
\def\Ints{\mbox{\rm Int}^s}

\def\Bmt{\overline{B_{\ep_0}(\mt)}}
\def\Lye{L_{y,\eta}}
\def\Lyep{L^{(p)}_{y,\eta}}
\def\Fyp{F^{(p)}_y}
\def\Fxp{F^{(p)}_x}
\def\Lxx{L_{x,\xi}}
\def\Lxxp{L^{(p)}_{x,\xi}}

\def\Wuo{W^{u,1}}
\def\Wui{W^{u,i}}
\def\Wuj{W^{u,j}}
\def\Wut{\tW^{u,2}}
\def\Wuk{W^{u,k}}
\def\Wuh{\hW^{u}}
\def\tWuo{\tW^{u,1}}
\def\tWui{\tW^{u,i}}
\def\tWuj{\tW^{u,j}}
\def\tWuk{\tW^{u,k}}
\def\hWuo{\hW^{u,1}}
\def\hWui{\hW^{u,i}}
\def\hWuj{\hW^{u,j}}
\def\hWuk{\hW^{u,k}}
\def\dj{\delta^{(j)}}
\def\do{\delta^{(1)}}
\def\epj{\ep^{(j)}}
\def\epo{\ep^{(1)}}
\def\hSj{\widehat{S}^{(j)}}
\def\hSo{\widehat{S}^{(1)}}

\def\tmu{\tilde{\mu}}
\def\tla{\tilde{\lambda}}
\def\hE{\widehat{E}}
\def\uk{u^{(k)}}
\def\ui{u^{(i)}}
\def\uj{u^{(j)}}
\def\vk{v^{(k)}}
\def\vl{v^{(l)}}
\def\vi{v^{(i)}}
\def\vj{v^{(j)}}
\def\wk{w^{(k)}}
\def\wi{w^{(i)}}
\def\wj{w^{(j)}}
\def\etak{\eta^{(k)}}
\def\etai{\eta^{(i)}}
\def\etaj{\eta^{(j)}}
\def\zetak{\zeta^{(k)}}
\def\zetai{\zeta^{(i)}}
\def\zetaj{\zeta^{(j)}}


\def\yj{y^{(j)}}
\def\yi{y^{(i)}}
\def\tyi{\ty^{(i)}}
\def\yo{y^{(1)}}
\def\zj{z^{(j)}}
\def\zo{z^{(1)}}
\def\vj{v^{(j)}}
\def\vo{v^{(1)}}
\def\kaj{\kappa^{(j)}}
\def\kao{\kappa^{(1)}}

\def\tyj{\tilde{y}^{(j)}}
\def\yl{y^{(\ell)}}
\def\tyl{\tilde{y}^{(l)}}
\def\wo{w^{(1)}}
\def\vo{v^{(1)}}
\def\vi{v^{(i)}}
\def\vj{v^{(j)}}
\def\vk{v^{(k)}}
\def\uo{u^{(1)}}
\def\wt{w^{(2)}}
\def\xio{\xi^{(1)}}
\def\xit{\xi^{(2)}}
\def\xii{\xi^{(i)}}
\def\xij{\xi^{(j)}}
\def\hxio{\hxi^{(1)}}
\def\hxit{\hxi^{(2)}}
\def\hxii{\hxi^{(i)}}
\def\hxij{\hxi^{(j)}}

\def\cxi{\check{\xi}}
\def\cxio{\cxi^{(1)}}
\def\cxit{\cxi^{(2)}}
\def\cet{\check{\eta}}
\def\ceto{\cet^{(1)}}
\def\cett{\cet^{(2)}}
\def\cv{\check{v}}
\def\cvo{\cv^{(1)}}
\def\cvt{\cv^{(2)}}
\def\cu{\check{u}}
\def\cuo{\cu^{(1)}}
\def\cut{\cu^{(2)}}
\def\cj{c^{(j)}}
\def\fj{f^{(j)}}
\def\gji{g^{(j,i)}}
\def\tPsi{\widetilde{\Psi}}
\def\chU{\check{U}}


\def\Ulo{U^{(\ell_0)}}
\def\dte{D_\theta}
\def\diamte{\mbox{\rm diam}_{\theta}}
\def\Ial{I^{(\alpha)}}
\def\uml{u_m^{(\ell)}}
\def\yl{y^{(\ell)}}
\def\tyl{\tilde{y}^{(\ell)}}
\def\ool{\oo^{(\ell)}}
\def\fl{f^{(\ell)}}
\def\hep{\hat{\ep}}
\def\dl{d^{(\ell)}}
\def\Lipt{{\Lip_\theta}}
\def\lipt{{\footnotesize \Lip_\theta}}
\def\tm{\tilde{m}}
\def\tj{\tilde{j}}
\def\lengthf{\mbox{\rm\footnotesize length}}
\def\length{\mbox{\rm length}}


\def\Xijl{X^{(\ell)}_{i,j}}
\def\hXijl{\widehat{X}^{(\ell)}_{i,j}}
\def\Wl{W^{(\ell)}}
\def\omijl{\omega^{(\ell)}_{i,j}}

\def\Xitl{X^{(\ell)}_{i,t}}
\def\hXitl{\widehat{X}^{(\ell)}_{i,t}}
\def\Vl{V^{(\ell)}}
\def\omitl{\omega^{(\ell)}_{i,t}}
\def\Xisl{X^{(\ell)}_{i,s}}
\def\hXisl{\widehat{X}^{(\ell)}_{i,s}}
\def\omisl{\omega^{(\ell)}_{i,s}}
\def\hGa{\widehat{\Gamma}}
\def\hOm{\widehat{\Omega}}
\def\tGa{\widetilde{\Gamma}}
\def\hA{\widehat{A}}
\def\tnu{\tilde{\nu}}
\def\tX{\widetilde{X}}

\def\ww{{\mathcal W}}
\def\Zl{Z^{(\ell)}}
\def\hpp{\widehat{\pp}}
\def\tnn{\widetilde{\nn}}

\def\fa{F^{(a)}}
\def\f0{F^{(0)}}
\def\tu{\tilde{u}}
\def\tD{\widetilde{D}}
\def\tchi{\tilde{\chi}}
\def\tC{\widetilde{C}}
\def\hC{\widehat{C}}
\def\hQ{\widehat{Q}}
\def\hF{\widehat{F}}
\def\hD{\widehat{D}}
\def\hr{\hat{r}}
\def\psid{\psi^\dag}
\def\taud{\tau^\dag}
\def\Omn{\Omega^{(n)}}
\def\Omm{\Omega^{(m)}}
\def\Omk{\Omega^{(k)}}
\def\Conf{{\mbox{\footnotesize\rm Const}}}
\def\hp{\hat{p}}

\def\tj{t^{(j)}}
\def\tyj{\tilde{y}^{(j)}}
\def\tyjo{\tilde{y}_{j,1}}
\def\tyjt{\tilde{y}_{j,2}}
\def\tyji{\tilde{y}_{j,i}}
\def\yjo{y_{j,1}}
\def\yjt{y_{j,2}}
\def\yji{y_{j,i}}
\def\tylo{\tilde{y}_{\ell,1}}
\def\tylt{\tilde{y}_{\ell,2}}
\def\tyli{\tilde{y}_{\ell,i}}
\def\ylo{y_{\ell,1}}
\def\ylt{y_{\ell,2}}
\def\yli{y_{\ell,i}}

\def\ulo{u_{\ell,1}}
\def\ult{u_{\ell,2}}
\def\uli{u_{\ell,i}}
\def\tulo{\tilde{u}_{\ell,1}}
\def\tult{\tilde{u}_{\ell,2}}
\def\tuli{\tilde{u}_{\ell,i}}

\def\vlo{v^{(\ell)}_{1}}
\def\vlt{v^{(\ell)}_2}
\def\vli{v^{(\ell)}_{i}}
\def\tvlo{\tilde{v}^{(\ell)}_{1}}
\def\tvlt{\tilde{v}^{(\ell)}_{2}}
\def\tvli{\tilde{v}^{(\ell)}_{i}}

\def\tdlo{\tilde{d}_{\ell,1}}
\def\tdlt{\tilde{d}_{\ell,2}}
\def\tdli{\tilde{d}_{\ell,i}}
\def\dlo{d_{\ell,1}}
\def\dlt{d_{\ell,2}}
\def\dli{d_{\ell,i}}
\def\wjo{w_{j,1}}
\def\wjt{w_{j,2}}
\def\wji{w_{j,i}}
\def\sj{s^{(j)}}
\def\Yj{Y^{(j)}}
\def\Vj{V^{(j)}}
\def\Zj{Z^{(j)}}
\def\vj{v^{(j)}}
\def\wj{w^{(j)}}
\def\twj{\tilde{w}^{(j)}}
\def\gj{g^{(j)}}
\def\tgj{\tilde{g}^{(j)}}
\def\tg{\tilde{g}}
\def\hn{\hat{n}}

\def\hbeta{\hat{\beta}}
\def\hmu{\hat{\mu}}
\def\piS{\pi^{(S)}}
\def\hb{\hat{b}}
\def\shP{P^\sharp}
\def\tshP{\widetilde{P}^\sharp}
\def\T{\mathcal T}
\def\tut{\tu^{(2)}}
\def\twt{\tw^{(2)}}
\def\piS{\pi^{(S)}}
\def\tsigma{\tilde{\sigma}}
\def\td{\tilde{d}}
\def\m{{\sf m}}
\def\Dom{D^{(\omega)}}
\def\Domo{D^{\omega_1)}}
\def\Qom{Q^{(\omega)}}
\def\Qomo{Q^{\omega_1)}}
\def\Aom{A^{(\omega)}}
\def\Aomo{A^{\omega_1)}}
\def\Aep{A^{(\ep)}}
\def\Aepp{A^{(\ep/2)}}
\def\W{{\mathcal W}}
\def\Ei{E^{(i)}}
\def\tEi{\widetilde{E}^{(i)}}

\begin{center}
{\Large\bf Lyapunov exponents and strong exponential tails\\ for some contact Anosov flows}
\end{center}

\begin{center}
{\sc by Luchezar Stoyanov}
\end{center}

\footnotesize

\noindent
{\bf Abstract.} For the time-one map $f$ of a contact Anosov flow on a compact Riemann manifold $M$,
satisfying a certain regularity condition,
we show that given a Gibbs measure on $M$, a sufficiently large Pesin regular set $P_0$ and an arbitrary
$\delta \in (0,1)$, there exist positive constants $C$ and $c$ such that for any integer $n \geq 1$, 
the measure of the set of those $x\in M$ with $f^k(x) \notin P_0$ for at least $\delta n$ values of $k = 0,1, \ldots,n-1$ 
does not exceed $C e^{-cn}$.

\normalsize

\section{Introduction}
\renewcommand{\theequation}{\arabic{section}.\arabic{equation}}

Let $\phi_t : M \longrightarrow M$ be a $C^2$ Anosov flow on a $C^2$ compact Riemann manifold 
$M$, and let $f = \phi_1$ be its {\it time-one map}.

It follows from a well-know result of Oseledets (\cite{kn:Os}; see also \cite{kn:BP} or
\cite{kn:R}) that there exists a Borel subset $\ll_0$ of $M$, which has full measure with
respect to any $f$-invariant Borel probability measure on $M$, such that for every 
$x\in \ll_0$ there exists a $df$-invariant decomposition
$$T_xM = E_1(x) \oplus E_2(x) \oplus \ldots \oplus E_{k(x)}(x)$$
and numbers $\chi_1 (x) < \chi_2(x) < \ldots < \chi_{k(x)}(x)$, called {\it Lyapunov characteristic exponents},
such that:

(a) $\di \lim_{|n| \to \infty} \frac{1}{n} \log \| df_x^n (v)\| = \chi_i(x)$ for all
$v \in E_i(x)\setminus \{0\}$ and all $i = 1, \ldots,k(x)$.

(b) For every $\ep > 0$ there exists a Borel function 
$A_\ep : \ll_0 \longrightarrow (1,\infty)$, such that
\be
\frac{\|v\|}{A_\ep(x)\, e^{|n|\ep}} \leq \frac{\|df_x^n(v)\|}{e^{n\chi_i(x)}} 
\leq A _\ep(x)\, e^{|n|\ep} \|v \| \quad , \quad  v\in E_i(x) \;, \; n \in \Z ,
\ee
for all $x\in \ll_0$ and all $i = 1, \ldots,k(x)$, and
\be
e^{-\ep} \leq \frac{A_\ep (f(x))}{A_\ep (x)} \leq  e^{\ep} \quad , \quad x\in \ll_0 .
\ee

(c) For all $x\in \ll_0$ and all disjoint non-empty subsets $I, I'$ of $\{1, \ldots, k(x)\}$ the smallest angle between non-zero
vectors in $E_I(x) = \oplus_{i\in I} E_i(x)$ and $E_{I'}(x)$ is $\geq \frac{1}{A_\ep(x)}$.

(d) If $\m$ is an ergodic $f$-invariant Borel probability measure on $M$, then the functions $k(x)$
and $\chi_i(x)$ are constant $\m$-almost everywhere.

\bs

A function $A_\ep$ satisfying (1.2) is called an {\it $\ep$-slow varying function}.

Let $\Phi$ be a H\"older continuous real-valued function on $M$ and let $\m$ be the {\it Gibbs 
measure} generated by $\Phi$ on $M$ (see e.g. \cite{kn:P2} or \cite{kn:PP}). Then $\m$ is ergodic, so 
there exists a subset $\ll'_0$ of $\ll_0$ with $\m(\ll'_0) = 1$ such that the functions $k(x) = k_0$ and
$\chi_i(x) = \chi_i$ are constant on $\ll'_0$.

It follows e.g. from the arguments in Sect. 3 in \cite{kn:PS}, that for every $\ep> 0$ 
there exist constants $\ep'\in (0,\ep]$, $p \geq 1$ and $\nu > 0$ such that 
$r_{\ep}(x) = \frac{\nu}{(A_{\ep'}(x))^p} $
determines an $\ep$-slowly varying function on $\ll'_0$ which defines a Lyapunov regular neighbourhood for 
every $x\in \ll'_0$, i.e. for each $x\in \ll'_0$ there exists a Lyapunov chart on $B(x,r_\ep(x))$. On these charts
one has estimates of the iterations of the non-linear map $f$ similar to these in (1.1).

It is known that in general the complement of the set $\ll_0$ (and therefore that of $\ll'_0$) can be 
topologically very large (see \cite{kn:BSau}, \cite{kn:BS} or \cite{kn:PSa} for some interesting 
examples). The regularity functions $A_\ep$ and $r_\ep$ are in general only measurable. 
The so called {\it Pesin regular sets}
$$\rr_\ell = \{ x\in \ll'_0 : A_\ep(x) \leq \ell\; , \; r_\ep(x) \geq 1/\ell\}  \quad , \quad \ell \geq 1 ,$$
and their closures are of particular importance since on such sets uniform estimates involving Lyapunov exponents
are available (see e.g. \cite{kn:P1}, \cite{kn:BP}, \cite{kn:KM}, \cite{kn:LY1}, \cite{kn:LY2}, \cite{kn:PS},
\cite{kn:BPS}). However it seems there is little information in the literature about the measures of the sets $\rr_\ell$ 
and `how quickly' they fill in $\ll'_0$, even in the case of uniformly hyperbolic systems.

Let $\m$ be an ergodic $f$-invariant Borel probability measure on $M$, and let $k_0$ be so that 
$k(x) = k_0$ for $\m$-almost all $x$.  Consider the distributions
$$\Ei(x) = E_1(x) \oplus E_2(x) \oplus \ldots \oplus E_i(x) \quad, \quad 
\tEi(x) = E_{i+1}(x) \oplus \ldots \oplus E_{k_0}(x) .$$
We will say that $\Ei(x)$ is
{\it uniformly continuous in $\ll'_0$}  if  the map $\ll'_0 \ni x \mapsto \Ei(x)$ is uniformly continuous with
respect to the natural distance between distributions with the same dimension (see e.g. Sect. 2.3 in \cite{kn:BP}).

Given a Lyapunov regularity function $A_\ep$, an integer $p \geq 0$ and a constant $\delta \in (0,1)$ set
$$T_p = \{ x\in \ll'_0 : A_\ep(x) \leq e^p\}, $$
and denote by $\Gamma_n = \Gamma_n(\ep,\delta,p)$ {\it the set of all $x\in \ll'_0$ with $f^k(x) \notin T_p$ for at least}
$\delta \, n$ values of $k = 0,1, \ldots,n-1$. 

\bs

In this paper we prove the following.

\bs

\noindent
{\bf Theorem 1.1.} {\it  Let $\phi_t : M \longrightarrow M$ be a $C^2$ contact  Anosov flow on a $C^2$ compact  
Riemann manifold $M$, let $\Phi$ be a H\"older continuous real-valued function on 
$M$ and let $\m$ be the {\it Gibbs measure} generated by $\Phi$ on $M$. Assume in addition that 
for every $i = 1, \ldots,k_0-1$ the distributions $\Ei(x)$ and $\tEi(x)$ are uniformly continuous in $\ll'_0$. Then
for every $\ep > 0$  there exists a Lyapunov $\ep$-regularity function $A_\ep$ satisfying {\rm (1.1)} and {\rm (1.2)} 
such that the following are satisfied:}

\ms

(a) (Exponential Tails) {\it There exist constants $C = C(\ep) > 0$, $c = c(\ep) > 0$ and $p_0 = p_0(\ep) > 0$ such that
$$\di \m\left(M \setminus \cup_{k = 0}^{n-1} f^{-k}(T_{p_0})  \right) \leq C\, e^{-cn} $$
for every integer  $n \geq 0$.}

\ms

(b) (Strong Exponential Tails) {\it For every $\delta \in (0,1)$ there exist constants 
$p_0 = p_0(\ep,\delta) > 0$, $C = C(\ep,\delta) > 0$ and  $c = c(\ep,\delta) > 0$ such that
$$\di \m\left(\Gamma_n(\ep,\delta,p_0)  \right) \leq C\, e^{-c n}$$
for every integer $n \geq 0$.}

\bs

 Clearly, taking the constant $p_0$ sufficiently large, we can make  $\m(T_{p_0})$ arbitrarily close to $1$.
The assumption that the flow is contact is used essentially in Sect. 6 below.

As we mentioned earlier, for every $p \geq 0$ there exists a constant $r_0 > 0$ (depending on $\ep$ and $p$) 
such that $r_\ep(x) \geq r_0$
for every $x \in T_p$. That is, for every $x\in T_p$ there exists a Lyapunov chart on $B(x,r_0)$.

The distributions $\Ei(x)$ and $\tEi(x)$ are uniformly continuous (in fact, H\"older continuous) 
for example when the flow has only three Lyapunov exponents
$\chi_1 < \chi_2 = 0 < \chi_3$, or more generally, when 
for every $i = 1, \ldots,  k_0-1$ there exist constants  $\lambda_i$, $\mu_i$  and $C > 0$ such that 
$\chi_i < \lambda_i  < \mu_i < \chi_{i+1}$,
$\| d\phi_t(u)\| \leq C\, \lambda_i^t\, \|u\|$ for all  $u\in \Ei(x)$ and $t\geq 0$, and
$\| d\phi_t(u)\| \geq \frac{1}{C}\, \mu_i^{t}\, \|u\|$ for all $u\in \tEi(x)$ and  $t\geq 0$
(see e.g. Ch. 3  in \cite{kn:P3}).

Thus, the conclusions in Theorem 1.1 always hold for contact Anosov flows on $3$-dimensional manifolds. As mentioned above,
a more general immediate consequence of Theorem 1.1 is the following.

\bs

\noindent
{\bf Corollary 1.2.} {\it  Let $\phi_t : M \longrightarrow M$ be a $C^2$ contact  Anosov flow on a $C^2$ compact  
Riemann manifold $M$, let $\Phi$ be a H\"older continuous real-valued function on 
$M$ and let $\m$ be the {\it Gibbs measure} generated by $\Phi$ on $M$. Assume that $k(x) = 3$ for almost all $x\in \ll_0$.
Then for every $\ep > 0$  there exists a Lyapunov $\ep$-regularity function $A_\ep$ satisfying {\rm (1.1)} and {\rm (1.2)} 
such that both {\rm (a)} and {\rm (b)} in Theorem 1.1  are satisfied.}

\bs

In this paper we restrict ourselves to contact Anosov flows, however using slight modifications of
the arguments in Sects. 3-5 below, results similar to Theorem 1.1 above can be proved for contact Anosov
diffeomorphisms\footnote{That case is in fact easier.}.

The motivation for the above result came from investigations on spectral properties of the so called 
Ruelle transfer operators. Attempts to obtain exponentially small estimates of integrals involving iterations
of a certain kind of `contraction operators' (see \cite{kn:D}, \cite{kn:St}) naturally lead to estimates
of `tails' of the kind considered in Theorem 1.1. On the other hand, studies on decay of correlations\footnote{See \cite{kn:L},
\cite{kn:BaL} and the references there for general information on this topic.} using the so called Young towers (\cite{kn:Y2})
also involve assumptions on the measures of the tails, and usually exponential tails are associated with
exponential decay of correlations. Such assumptions appear naturally also in studies on large deviations (see e.g.
\cite{kn:MN} and the references there).






Sect. 2 below contains some basic definitions. In Sect. 3 we introduce a `regularity function' $R_\ep(x)$ 
(and $\tR_\ep(x)$ related to it),
changing slightly the approach of Simic in \cite{kn:Sim}. Although this is most likely not a slowly varying function,
it turns out that $\con(\sqrt{\ep} A_{2\ep}(x))^{1/3} \leq R_\ep(x) \leq A_\ep(x)$ for some
Lyapunov $\ep$-regularity function  $A_\ep(x)$ satisfying
(1.1) and (1.2). This is proved in Sect. 6. Consequently, it is enough to prove the analogue
of Theorem 1.1 replacing $A_\ep(x)$ by $R_\ep(x)$. We study the latter (or rather, the related one $\tR_\ep(x)$)
in Sects. 4 and 5 using an idea\footnote{We use the idea briefly mentioned at the end of Sect. 1 in \cite{kn:Sim}
apparently suggested by the referee of that paper.} from \cite{kn:Sim} and the classical large deviation principle 
(see \cite{kn:OP}, \cite{kn:Ki} or \cite{kn:Y1}). Theorem 1.1 is proved in Sect. 4 as a consequence of a
similar result about the function $\tR_\ep(x)$.

\bs

\footnotesize
\noindent
{\bf Acknowledgement.} Thanks are due to Sebastian Gou\"ezel who pointed out to an error in a
previous version of the paper.

\normalsize


\def\ho{h^{(1)}}
\def\Ho{H^{(1)}}

\section{Preliminaries}
\setcounter{equation}{0}

Throughout this paper $M$ denotes a $C^2$ compact Riemann manifold,  and 
$\phi_t : M \longrightarrow M$ ($t\in \R$) a  $C^2$ Anosov flow on $M$. That is, for some
constants $C > 0$ and $0 < \lambda < 1$ there exists a $d\phi_t$-invariant decomposition  
$T_xM = E^0(x) \oplus E^u(x) \oplus E^s(x)$ of $T_xM$ ($x \in M$) 
into a direct  sum of non-zero linear subspaces,
where $E^0(x)$ is the one-dimensional subspace determined by the direction of the flow
at $x$, $\| d\phi_t(u)\| \leq C\, \lambda^t\, \|u\|$ for all  $u\in E^s(x)$ and $t\geq 0$, and
$\| d\phi_t(u)\| \leq C\, \lambda^{-t}\, \|u\|$ for all $u\in E^u(x)$ and  $t\leq 0$.
For $x\in M$ and a sufficiently small $\epsilon > 0$ let 
$$\wloc^s(x) = \{ y\in M : d (\phi_t(x),\phi_t(y)) \leq \epsilon \: \mbox{\rm for all }
\: t \geq 0 \; , \: d (\phi_t(x),\phi_t(y)) \to_{t\to \infty} 0\: \}\; ,$$
$$\wloc^u(x) = \{ y\in M : d (\phi_t(x),\phi_t(y)) \leq \epsilon \: \mbox{\rm for all }
\: t \leq 0 \; , \: d (\phi_t(x),\phi_t(y)) \to_{t\to -\infty} 0\: \}$$
be the (strong) {\it stable} and {\it unstable manifolds} of size $\epsilon$. Then
$E^u(x) = T_x \wloc^u(x)$ and $E^s(x) = T_x \wloc^s(x)$.

The flow $\phi_t$ is called  {\it contact} if $\dim(M) = 2n+1$
for some $n \geq 1$ and there exists a $C^2$ flow-invariant one-form $\omega$ on $M$ 
such that $\omega \wedge (d\omega)^n \neq 0$ on $M$. It is well-known that the 
Lyapunov spectrum of a contact flow is symmetric, i.e. for each $i = 1, \ldots, k_0$ 
there exists $j = 1, \ldots, k_0$  with $\chi_j = - \chi_i$. 

It follows from the hyperbolicity of the flow on $M$  that if  $\epsilon_0 > 0$ is sufficiently small,
there exists $\ep > 0$ such that if $x,y\in M$ and $d (x,y) < \ep$, 
then $W^s_{\ep_0}(x)$ and $\phi_{[-\ep_0,\ep_0]}(W^u_{\ep_0}(y))$ intersect at exactly 
one point $[x,y]$  (cf. \cite{kn:KH}). That is, there exists a unique  $t\in [-\ep_0, \ep_0]$ such that
$\phi_t([x,y]) \in W^u_{\ep_0}(y)$.

Let $\rr = \{ R_i\}_{i=1}^{k_0}$ be a Markov family consisting of rectangles $R_i$, each contained
in a submanifold $D_i$ of $M$ of codimension one (see \cite{kn:B}). 
Assuming that $\ep > 0$ is sufficiently small, the projection 
$\pr_{D_i} : \phi_{[-\ep,\ep]}(D_i) \longrightarrow D_i$ along the flow is well-defined
and smooth. Given $x,y\in D_i$, set $\la x, y\ra_{D_i} = \pr_{D_i}([x,y])$. 
A subset $R_i$ of $D_i$ is called a {\it rectangle} if $\la x, y\ra_{D_i} \in R_i$ for all $x,y\in R_i$.  
For any $x\in R_i$ define the stable and unstable leaves through $x$ in $R_i$
by $W^s_{R_i}(x) = \pr_{D_i}(W^s_\ep(x)\cap \phi_{[-\ep,\ep]}(D_i)) \cap R_i$ and
$W^u_{R_i}(x) = \pr_{D_i}(W^u_\ep(x)\cap \phi_{[-\ep,\ep]}(D_i)) \cap R_i$. 
We may assume that each $R_i$ has the form 
$R_i = \la U_i  , S_i \ra_{D_i} = \{ \la x,y\ra_{D_i} : x\in U_i, y\in S_i\}$,
where $U_i \subset \wloc^u(z_i)$ and $S_i \subset \wloc^s(z_i)$, respectively,  for some $z_i\in M$.  
Set $R =  \cup_{i=1}^{k_0} R_i$. 
The corresponding {\it Poincar\'e map} $\pp: R \longrightarrow R$ is defined by  
$\pp(x) = \phi_{\tau(x)}(x) \in R$, where $\tau(x) > 0$ is the smallest positive time with $\phi_{\tau(x)}(x) \in R$. 
The function $\tau$  is the {\it first return time}  associated with $\rr$. 

From now on we will assume that $\rr = \{ R_i\}_{i=1}^{k_0}$ is a fixed Markov family for  
$\phi_t$ of size $\chi < \ep_0/4 < 1/4$. 
Denote by $\hR$ the {\it core} of  $R$, i.e. the set $x \in R$ such that  $\pp^m(x) \in \Int(R) = \cup_{i=1}^k \Int(R_i)$  
for all $m \in \Z$. It is well-known (see \cite{kn:B}) that $\hR$ is a residual subset 
of $R$ and has full measure with respect to any Gibbs measure on $R$.
In general $\tau$ is not continuous on $R$, however it is H\"older on $\hR$ when considered with respect to
an appropriate metric $d_\theta$, see below. The same applies to $\pp : \hR \longrightarrow \hR$.  
It is well-known (\cite{kn:BR}) that the Markov family $\rr$ can be chosen so that
{\bf $\tau$ is non-lattice}. From now on we will assume that $\rr$ is chosen in this way.

Let $\aa = (\aa_{ij})_{i,j=1}^k$ be the matrix given by $\aa_{ij} = 1$ 
if $\pp(\Int (R_i)) \cap \Int (R_j) \neq  \e$ and $\aa_{ij} = 0$ otherwise. 
Consider the  symbol space 
$$\sa = \{  (i_j)_{j=-\infty}^\infty : 1\leq i_j \leq k_0, A_{i_j\; i_{j+1}} = 1
\:\: \mbox{ \rm for all } \: j\; \},$$
with the product topology and the {\it shift map} $\sigma : \sa \longrightarrow \sa$ 
given by $\sigma ( (i_j)) = ( (i'_j))$, where $i'_j = i_{j+1}$ for all $j$.
Given $0 < \theta < 1$, consider  the {\it metric} $d_\theta$ on $\sa$ defined by
$d_\theta(\xi,\eta) = 0$ if $\xi = \eta$ and $d_\theta(\xi,\eta) = \theta^m$ if
$\xi_i = \eta_i$ for $|i| < m$ and $m$ is maximal with this property. 
There is a  natural map $\W : \sa \longrightarrow R$ such that
$\W\circ \sigma = \pp\circ \W$. In general $\W$ is not one-to-one, however it is a bijection between
$\W^{-1}(\hR)$ and $\hR$.  Choosing $\theta \in (0,1)$ appropriately, the map 
$\W: \sa \longrightarrow R$ is Lipschitz when $\sa$ is endowed with $d_\theta$ and $R$ with
the metric induced by the Riemann metric on $M$. In what follows we assume $\theta \in (0,1)$
is a fixed  constant with this property, and we will consider $\hR$ with the metric induced by
$d_\theta$ via $\W$.

\section{Lyapunov regularity functions}
\setcounter{equation}{0}

Throughout we assume that $\phi_t$ is a $C^2$ contact Anosov flow on $M$. 
Let $\Phi : M \longrightarrow \R$ be a fixed H\"older continuous function on $M$  and let
$\m$ be the {\it Gibbs measure} determined by $\Phi$. Let $\ll'_0$ be a {\it subset of $\ll_0$ of full $\m$-measure} 
such that $k(x) = k_0$ and $\chi_i(x)$ are constant for $x\in \ll'_0$ for all $i = 1, \ldots, k_0$.
Removing a set of measure zero, we may assume that $\ll'_0 \subset \hR$.

Take  $\ep > 0$ so small that $\chi_j \notin [\chi_i - 4\ep, \chi_i + 4\ep]$ whenever $i \neq j$, and set
\be
R^+_{\ep}(x) = \max_{1\leq i \leq k_0} \: \sup_{n\geq 0} \; \frac{\| (df^n_x)_{| E_i}\|}{e^{(\chi_i + \ep) n}} ,
\ee
\be
R^-_{\ep}(x) = \max_{1\leq i \leq k_0} \: \sup_{n\geq 0} \; \frac{\| (df^{-n}_x)_{| E_i}\|}{e^{(\chi_i - \ep) (-n)}} ,
\ee
and
$R_\ep(x) = \max \{ R^-_\ep(x), R^+_\ep(x) \}$.

\bs

\noindent
{\bf Note.} Since $\frac{1}{n} \log \| (df^n_x)_{| E_i}\| \to \chi_i$ as $n \to \infty$, we have
$\frac{1}{n} \log \| (df^n_x)_{| E_i}\| <  \chi_i + \ep$ for large $n$, so (3.1) exists.
Similarly,  $\frac{1}{-n} \log \| (df^{-n}_x)_{| E_i}\| \to \chi_i$ as $n \to \infty$, so
$\frac{1}{-n} \log \| (df^{-n}_x)_{| E_i}\| >  \chi_i - \ep$ for large $n$.
This gives  $\log \| (df^{-n}_x)_{| E_i}\| < (\chi_i-\ep) (-n)$ for large $n$, and so  (3.2) exists.

\bs

From the above definitions it is clear that 
\be
\frac{\|df_x^n(v)\|}{e^{n\chi_i} } \leq R _\ep(x)\, e^{|n|\ep} \|v\|\quad , \quad v\in E_i(x)
\;, \; n \in \Z ,
\ee
for all $x\in \ll'_0$ and all $i =1, \ldots,k_0$, 
and moreover for every Lyapunov $\ep$-regularity function $A_\ep(x)$ with (1.1) we have
$R_\ep(x) \leq A_\ep(x)$
for all $x\in \ll'_0$ and $\ep > 0$. It is not clear whether $R_\ep(x)$ is a slowly varying function, i.e.
whether it satisfies (1.2), possibly with $\ep$ replaced by some $\delta > 0$ related to $\ep$. 
However we have the following.

\bs

\noindent
{\bf Proposition 3.1.} {\it  Let $\phi_t : M \longrightarrow M$ be a contact $C^2$ Anosov flow,
and let $\Phi$, $\m$ and $\ll'_0$ be as above. For every $\ep > 0$ there exists an  $\ep$-slowly varying 
function  $A_\ep$ satisfying {\rm (1.1)} and {\rm (1.2)} such that 
\be
\frac{\ep^{1/6}}{C} (A_{2\ep}(x))^{1/3} \leq R_\ep(x) \leq A_\ep(x) 
\ee
for almost all $x\in \ll'_0$, where  $C = C_0^{2/3} (2k_0)^{1/6} > 1$.} 

\bs

Notice that $ \frac{\ep^{1/6}}{C} (A_{2\ep}(x))^{1/3}$ is also an $\ep$-slowly varying function.

We prove this proposition in Sect. 6 below. The assumption that the flow is contact
is used in an essential way there.

We will also need the functions
$$\tR^+_{\ep}(x) = \max_{1\leq i \leq k_0} \: \sup_{n\geq 0} \; \frac{\| (df^n_x)_{| \Ei}\|}{e^{(\chi_i + \ep) n}}  \quad , \quad
\tR^-_{\ep}(x) = \max_{1\leq i \leq k_0} \: \sup_{n\geq 0} \; \frac{\| (df^{-n}_x)_{| \tEi}\|}{e^{(\chi_i - \ep) (-n)}} ,$$
and
$\tR_\ep(x) = \max \{ \tR^-_\ep(x), \tR^+_\ep(x) \}$. Clearly $R_\ep(x) \leq \tR_\ep(x)$. In Sect. 6 below we will
show that there exists a constant $K(\ep) \geq 1$ such that 
\be
\tR_\ep(x) \leq K(\ep) \, (R_\ep(x))^2
\ee
for all $x \in \ll'_0$.

\bs

\noindent
{\bf Remark 3.2.}
Proposition 3.1 and (3.5) show that sets of the form $\{ x\in \ll'_0 :  A_\ep(x) \leq e^p\}$ are easily related to
sets of the form $\{ x\in \ll'_0 :  \tR_{\ep'}(x) \leq e^{p'}\}$ (with $\ep' = \ep$ or $3\ep$), so it is enough to 
prove an analogue of Theorem 1.1 replacing the function $A_\ep$ by $ C \tR_\ep$ for some constant $C > 0$.

\def\tyj{\tilde{y}^{(j)}}
\def\yl{y^{(\ell)}}
\def\wo{w^{(1)}}
\def\vo{v^{(1)}}
\def\vi{v^{(i)}}
\def\vj{v^{(j)}}
\def\vk{v^{(k)}}
\def\uo{u^{(1)}}
\def\wt{w^{(2)}}
\def\xio{\xi^{(1)}}
\def\xit{\xi^{(2)}}
\def\xii{\xi^{(i)}}
\def\xij{\xi^{(j)}}
\def\hxio{\hxi^{(1)}}
\def\hxit{\hxi^{(2)}}
\def\hxii{\hxi^{(i)}}
\def\hxij{\hxi^{(j)}}

\def\cxi{\check{\xi}}
\def\cxio{\cxi^{(1)}}
\def\cxit{\cxi^{(2)}}
\def\cet{\check{\eta}}
\def\ceto{\cet^{(1)}}
\def\cett{\cet^{(2)}}
\def\cv{\check{v}}
\def\cvo{\cv^{(1)}}
\def\cvt{\cv^{(2)}}
\def\cu{\check{u}}
\def\cuo{\cu^{(1)}}
\def\cut{\cu^{(2)}}
\def\cj{c^{(j)}}
\def\fj{f^{(j)}}
\def\gji{g^{(j,i)}}

\def\ao{a^{(1)}}
\def\bo{b^{(1)}}
\def\zi{z^{(i)}}
\def\hGa{\widehat{\Gamma}}

\def\tvo{\tv^{(1)}}
\def\tetao{\teta^{(1)}}

\def\fa{F^{(a)}}
\def\f0{F^{(0)}}
\def\tu{\tilde{u}}
\def\tD{\widetilde{D}}
\def\tchi{\tilde{\chi}}
\def\tC{\widetilde{C}}
\def\hC{\widehat{C}}
\def\hF{\widehat{F}}
\def\hD{\widehat{D}}
\def\hr{\hat{r}}
\def\psid{\psi^\dag}
\def\taud{\tau^\dag}
\def\Omn{\Omega^{(n)}}
\def\Conf{{\mbox{\footnotesize\rm Const}}}
\def\hp{\hat{p}}
\def\Ri{R^{(i)}}

\def\Pi{P^{(i)}}
\def\Xxi{X^{(i)}}

\def\Yi{Y^{(i)}}
\def\Zi{Z^{(i)}}
\def\Gn{G^{(n)}}
\def\Gso{G^{(s_0)}}
\def\tY{\widetilde{Y}}
\def\tll{\widetilde{\ll}}

\section{Reductions}
\setcounter{equation}{0}

Let $\phi_t$ be a $C^2$ Anosov flow on $M$, and let $\Phi$, $\m$ and $\ll'_0 \subset \ll_0$ be as in Sect. 3.  
In this section and  also the next one  we do not assume that the flow is contact. However, as in Theorem 1.1, 
we will assume that the bundles $\Ei(x)$ and $\tEi(x)$ are uniformly continuous on $\ll'_0$. 

Let $\rr = \{ R_i\}_{i=1}^{k_0}$  be a Markov family for $\phi_t$ as in Sect. 2, and let
$\tau: R = \cup_{i=1}^{k_0} R_i \longrightarrow [0,1/4]$ and $\pp : R\longrightarrow R$ be the corresponding first return map
and the Poincar\'e map.  Fix constants $0 < \tau_0 < \htau_0 < 1/4$ so that 
$\tau_0 \leq \tau(x) \leq \htau_0$  for all $x\in R$. There is a well-defined {\it projection}
$\pi : M \longrightarrow R$ defined by $\pi(y) = x$, where $x\in R$, $y = \phi_s(x)$ for some $s\in [0,\tau(x))$
and $s \geq  0$ is the smallest number with this property.  

The Gibbs measure $\m$ induces a {\it Gibbs measure} $\mu$ on $R$ (with respect to the Poincar\'e map $\pp$) for the function
$$F(x) = \int_0^{\tau(x)} \Phi(\phi_s(x))\, ds \quad, \quad x\in R .$$
The function $F$ is H\"older on $\hR$ (with respect to the metric $d_\theta$).
For every continuous function $H$ on $M$ we have (see e.g. \cite{kn:PP})
\be
\int_M H\, d\m = \frac{\int_R \left( \int_0^{\tau(x)} H(\phi_s(x))\, ds\right) d\mu(x)}{\int_R \tau\, d\mu } .
\ee

By Birkhoff's Theorem,
\be
\lim_{n\to \infty} \frac{\tau_n(x)}{n} =  \taud = \int_R \tau\, d\mu 
\ee
for $\mu$-almost all $x\in \ll'_0 \cap R$. Here $\tau_n(x) = \tau(x) + \tau(\pp(x)) + \ldots + \tau(\pp^{n-1}(x))$.
By the choice of $\tau$ we have $\taud \leq 1/4$.

{\bf Fix a constant $C_1 \geq 1$ so that}
$$C_1 \geq \frac{k_0}{\taud} \sup \{ | \log \| d\phi_t(x)\| | : x\in M\;, \; 0\leq t\leq 1\} ,$$
and $C_1 \geq 2 \max_{1\leq i \leq k_0} |\chi_i|$. Then
{\bf fix arbitrary constants }
\be
0 < \ep \leq \frac{1}{8} \min \{ |\chi_i| : 1\leq i \leq k_0, \chi_i \neq 0 \}  \quad, \quad \delta \in (0,1) \quad, \quad
C_0 \geq \frac{100 C_1}{\ep \tau_0}  ,
\ee
and set
\be
\delta_0 = \frac{\delta}{2k_0} .
\ee

We will study in details  the regularity function $\tR^+_\ep(x)$, In a similar way one can deal with
$\tR^-_\ep(x)$; one just needs to replace $\phi_t$ by $\phi_{-t}$.
Clearly $\tR_\ep(x) = \max_{1\leq i \leq k_0} \tR^+_{\ep,i}(x)$, where 
$$\tR^+_{\ep,i}(x) = \sup_{n \geq 0} \frac{\|(d f^n_x)_{|\Ei}\|}{e^{(\chi_i+\ep)n}} .$$


In this section, and the next section as well, we will restrict our attention to an arbitrary {\bf fixed $i = 1, \ldots, k_0$, and the 
corresponding bundle $E = \Ei$}. {\bf Set $\chi = \chi_i$ and $B_\ep(x) = \tR^+_{\ep,i}(x)$} for brevity.
Thus,
\be
B_\ep(x) = \sup_{n \in \N} \frac{\| d \phi_n(x)_{| E}\|}{e^{(\chi + \ep) n}} \quad, \quad x\in \ll'_0 ,
\ee
where $\N$ is the set of non-negative integers. 
We will compare $B_\ep(x)$ with the function
$$C_\ep(x) = \sup_{k \geq 0} \frac{\| d \phi_{\tau_k(x)}(x)_{| E}\|}{e^{(\chi + \ep) \tau_k(x)}} \quad, \quad x\in \ll'_0\cap R  .$$

\ms

\noindent
{\bf Lemma 4.1.}  (a) {\it There exists a constant $S_1 = S_1(\ep)  > 0$ such that
\be
B_\ep(y) \leq S_1\, C_\ep(x)
\ee
for all $y \in \ll'_0$, $x = \pi(y)$. Moreover we can take $S_1 > 0$ so that {\rm (4.6)} holds for any $y\in \ll$
and $x\in R$ with $y = \phi_t(x)$ for some $t \in [-2,2]$.}

\ms

(b) {\it For any integer $s_0 \geq 1$ there exists a constant $S_2 = S_2(\ep, s_0) \geq 1$ such that
$C_\ep(\pp^{-j}(x)) \leq S_2 \, C_\ep(x)$ for any $x\in \ll'_0 \cap R$ and any $j = 0,1,\ldots, s_0$.}

\ms 

\noindent
{\it Proof.} Set $a_1 = 2 (|\chi| + \ep)$ and $a_2 = \sup_{|s|\leq 2} \| d\phi_s\| < \infty$.

(a) Let $y \in \ll'_0$; set $x = \pi(y) \in \ll'_0 \cap R$. Then $y = \phi_s(x)$ for some $s\in [0,\tau(x))$. We have
$B_\ep(y) = \frac{\| d \phi_t(x)_{| E}\|}{e^{(\chi + \ep) t}}$ for some integer $t \geq 0$.
Denote by $k \geq 0$ the maximal integer so that $\tau_k(x) \leq t+s$; then $t+s < \tau_{k+1}(x)$. 
Thus, $t \geq \tau_k(x) -s \geq \tau_k(x) - \htau_0$. Also, 
\begin{eqnarray*}
\|d\phi_t(y)_{|E}\| 
& =    & \|d\phi_{-s + (t+s)}(\phi_s(x))_{|E}\| 
\leq \|d\phi_{-s}(\phi_s(x))_{|E}\| \cdot  \|d\phi_{t+s}(x)_{|E}\|\\
& \leq & a_2 \, \|d\phi_{\tau_k(x) + (t+s-\tau_k(x))}(x)_{|E}\|
\leq  a^2_2  \|d\phi_{\tau_k(x)}(x)_{|E}\| .
\end{eqnarray*}
Thus, 
\begin{eqnarray*}
B_\ep(y)
 =     \frac{\|d\phi_t(y)_{|E}\|}{e^{(\chi+\ep) t}}
\leq \frac{a^2_2\, \|d\phi_{\tau_k(x)}(x)_{|E}\|}{e^{(\chi + \ep) (\tau_k(x) - \htau_0)}} 
 \leq  e^{a_1} a_2^2\, \frac{ \|d\phi_{\tau_k(x)}(x)_{|E}\|}{e^{(\chi + \ep) \tau_k(x)}} .
\end{eqnarray*}
This shows that $B_\ep(y) \leq \Con\, C_\ep(x)$. In a similar way one proves that $C_\ep(x) \leq \Con\, B_\ep(y)$.

The more general case when $y = \phi_t(x)$ for some $t \in [-2,2]$ follows similarly.

\ms

(b) Take $S_2 = \max_{0\leq j\leq s_0} \sup_{z\in R}  \frac{\| d \phi_{\tau_j(z)}(z)\|}{e^{(\chi + \ep) \tau_j(z)}}$.
Given $x\in \ll'_0\cap R$ and $j = 0,1,\ldots, s_0$, set $y = \pp^{-j}(x)$. For  any integer $n \in[0,s_0]$ clearly
$\frac{\| d \phi_{\tau_n(y)}(y)_{| E}\|}{e^{(\chi + \ep) \tau_n(y)}} \leq S_2 \leq S_2 C_\ep(x)$. Let $n > s_0$. Then
$n = j+k$ for some integer $k > 0$, and therefore
\begin{eqnarray*}
\frac{\| d \phi_{\tau_n(y)}(y)_{| E}\|}{e^{(\chi + \ep) \tau_n(y)}}
& =    & \frac{\| d \phi_{\tau_{j+k}(y)}(y)_{| E}\|}{e^{(\chi + \ep) \tau_{j+k}(y)}} =
\frac{\| d \phi_{\tau_{j}(y) + \tau_k(x)}(y)_{| E}\|}{e^{(\chi + \ep) (\tau_{j}(y) + \tau_k(x))}}\\
& \leq & \frac{\| d \phi_{\tau_{j}(y)}(y)_{| E}\|}{e^{(\chi + \ep) \tau_{j}(y)} }
\cdot \frac{\| d \phi_{\tau_k(x)}(x)_{| E}\|}{e^{(\chi + \ep) \tau_k(x) }} \leq S_2\, C_\ep(x) .
\end{eqnarray*} 
Thus, $C_\ep(y) \leq S_2\, C_\ep(x)$.
\endofproof

\bs

\noindent


Next, consider the functions
$$\Gn(x) = \log \| d\phi_{\tau_n(x)}(x)_{| E}\|$$ 
for $x\in \ll'_0 \cap R$ and any $n \geq 0$.
Clearly, $|G^{(1)}(x)| \leq \Con < \infty$ for all $x\in \ll'_0 \cap R$, and
\be
\lim_{n \to \infty} \frac{\Gn(x)}{\taud n} = \lim_{n \to \infty} \frac{\Gn(x)}{\tau_n(x)} \, \frac{\tau_n(x)}{\taud n} 
= \chi 
\ee
for  $\mu$-almost all $x \in R$. Here we used the fact that $\lim_{t\to \infty} \frac{\log \|d\phi_t(x)_{| E}\|}{t} = \chi$
for $\mu$-almost all $x\in R$.

The functions $\Gn$ form a sub-additive sequence with respect to the map 
$\pp : R \longrightarrow R$. Indeed, for any $n,m \geq 0$ and any $x \in \ll'_0 \cap R$ we have
\begin{eqnarray*}
G^{(n+m)}(x)
& =    & \log \| d\phi_{\tau_{n+m}(x)}(x)_{| E}\| = \log \| d\phi_{\tau_{n}(x) + \tau_m(\pp^n(x))}(x)_{| E}\|\\
& \leq & \log \| d\phi_{\tau_{n}(x)}(x)_{| E}\| + \log \| d\phi_{\tau_m(\pp^n(x))}(\pp^n(x))_{| E}\| 
 =  \Gn(x) + G^{(m)}(\pp^n(x)) .
\end{eqnarray*}
By Kingsman's Ergodic Theorem (\cite{kn:Kin}) there exists the limit
\be
L = \lim_{n\to \infty} \frac{1}{n} \int_R \Gn(x) \, d\mu(x) =
\inf_{n} \frac{1}{n} \int_R \Gn(x) \, d\mu(x) ,
\ee
and moreover $\lim_{n\to\infty} \frac{\Gn(x)}{n} = L$ for $\mu$-almost all $x\in \ll'_0 \cap R$. Now (4.7) shows that $L = \taud \chi$.

{\bf Fix a subset $\ll$ of $\ll'_0 \cap R$ with} $\mu(\ll) = 1$ such that (4.2) and (4.7) hold for all $x\in \ll$.
Set $\tll = \pi^{-1}(\ll) \subset \ll'_0$; then $\tll \cap R = \ll$.

By Egorov's Theorem, there exists a compact subset $K_0$ of $\ll$ such that
\be
\mu(K_0) > 1 - \frac{\delta_0}{C_0} ,
\ee
and $\frac{\tau_n(x)}{\taud n} \to 1$ and  $\frac{\Gn(x)}{\taud n} \to \chi$ uniformly on $K_0$. 
Thus, there exists an integer $n_0 \geq 1$ such that
\be
\left| \frac{\tau_n(x)}{\taud n} - 1 \right| < \frac{\delta_0}{C_0} \quad ,\quad x\in K_0\;, 
\ee
\be
\left| \frac{\Gn(x)}{\taud n} - \chi \right| < \frac{\delta_0}{C_0} \quad ,\quad x\in K_0\;, 
\ee
and, using (4.8),
\be
\chi \leq \frac{1}{\taud n} \int_R \Gn(x) \, d\mu(x) < \chi + \frac{\delta_0}{C_0} 
\ee
for all integers $n \geq n_0$.

{\bf Fix an arbitrary integer $s_0 = s_0(\ep, \delta_0) \geq n_0$}; this will stay fixed throughout the whole Sects. 4 and 5.
Then (4.10), (4.11) and (4.12) hold with $n$ replaced by $s_0$.
Consider the transformation 
$$T = \pp^{s_0}: R \longrightarrow R$$  
preserving the measure $\mu$, and the measurable functions
$$u(x) = \frac{\Gso (x)}{\taud s_0}  = \frac{1}{\taud s_0} \log \| d\phi_{\tau_{s_0}(x)}(x)_{|E}\| \quad , \quad x\in \ll ,$$
and 
$$g(x) = \frac{\tau_{s_0} (x)}{\taud s_0}  \quad , \quad x\in \ll .$$
Clearly $u$ and $g$ are bounded, $\|u\|_\infty \leq C_1$ and $\tau_0/\taud \leq g(x) \leq \htau_0/\taud$ for all $x\in R$.
For any integer $n \geq 1$ and $x\in \ll$ set
$$u^n(x) = u(x) + u(T(x)) + \ldots + u(T^{n-1}(x)) .$$
(Notice that a superscript is used here when we deal with orbits with respect to $T = \pp^{s_0}$; unlike the 
subscript used for orbits with respect to $\pp$.) 

It follows from our assumptions that $u$ is uniformly continuous on $\ll'_0$.
For the integral of $u$, (4.12) implies
\be
\chi \leq \int_R u \, d\mu < \chi + \frac{\delta_0}{C_0}  .
\ee



\noindent
{\bf Lemma 4.2.}  {\it There exist a Lipschitz\footnote{H\"older continuity is enough.} function $\psi : R \longrightarrow \R$ 
with respect to the metric $d_\theta$ such that $\|\psi\|_\infty \leq 2C_1$ and
$u(x) \leq \psi(x) \leq u(x) + \delta_0/C_0$ for $\mu$-almost all $x\in R$. }

\bs

\noindent
{\it Proof.} Since $u$ is uniformly continuous on $\ll'_0$ (with respect to the metric $d_\theta$), 
it has a continuous extension to the whole of $R$. We will denote this extension by $u$ again, and we
still have $\|u\|_\infty \leq C_1$. 

 Take $a = \frac{\delta_0}{2C_0}$ and consider the function $u + a$ on $R$. 
 On any given rectangle $R_i$ this function is continuous both with respect to the metric $d_\theta$ and the Riemann metric.
Using the manifold structure of the disk $D_i$ containing $R_i$ and a standard regularization procedure, we find a Lipschitz 
(with respect to the Riemann metric) function $\psi$ on $R_i$ such that $\|\psi - (u+ a)\|_\infty < a$ on $R_i$. 
Then $\psi$ is Lipschitz on $R$ with respect to $d_\theta$ and $\psi \geq u$. Moreover, $\|\psi - (u+ a)\|_\infty < a$ and 
$|u\|_\infty \leq C_1$ imply $\|\psi\|_\infty \leq 2C_1$. 
\endofproof

\bs

{\bf Fix a function} $\psi$ with the properties in Lemma 4.2. It follows from 
(4.13) that 
\be
\chi \leq  \tchi = \int_R \psi \, d\mu < \chi + \frac{2 \delta_0}{C_0} .
\ee
Changing slightly $\psi$ if necessary, we may assume that $\psi$ {\bf is not cohomologous to a constant}.

Next, set 
\be
D_\ep(x) =  \sup_{m\geq 0} \frac{e^{\psi^m(x)}}{e^{(\tchi + \ep) g^m(x)}} .
\ee
By Birkhoff's Theorem, this is well-defined for $\mu$-almost all $x\in \ll$. Removing a set of measure
zero from $\ll$  if necessary (and thus shrinking the compact set $K_0$ a bit), we will assume that $D_\ep(x)$
is defined for all $x \in \ll$. 

\ms

\noindent
{\bf Lemma 4.3.}  {\it There exists a constant $S_3 = S_3(\ep, \delta_0) > 0$ such that
$C_{2\ep}(x) \leq S_3 \, (D_\ep(x))^{\taud s_0}$ for all $x \in \ll$.}

\bs 

\noindent
{\it Proof.} Set $a_3 = \sup_{0 \leq t \leq s_0 \htau_0} \log \|d\phi_t(x)\|$.
Then $G^{(\ell)} (x) \leq a_3$ for all $0 \leq \ell \leq s_0$ and all $x\in \ll$.

Given $x\in \ll$, there exists $k \in \N$ such that 
$C_{2\ep}(x) = \frac{\|d\phi_{\tau_k(x)}(x)_{|E}\|}{e^{(\chi+2\ep) \tau_k(x)}}$.
Let $k = ms_0 + \ell$, $0 \leq \ell < s_0$. Then, setting $y = \pp^{ms_0}(x)$ and using $u \leq \psi$, we get
\begin{eqnarray*}
\log \|d\phi_{\tau_k(x)}(x)_{|E}\|
& =    & G^{(ms_0 + \ell)}(x) \leq \Gso(x) + \Gso(\pp^{s_0}(x)) + \ldots + \Gso(\pp^{(m-1)s_0}(x)) + G^{(\ell)}(y)\\
& \leq & \taud s_0 (u^m(x) + a_3) \leq \taud s_0 (\psi^m(x) + a_3) .
\end{eqnarray*}
Since $\tau_k(x) \geq \tau_{ms_0}(x)$  and $\tchi < \chi + \ep$ by (4.14) and the choice of the constants,
it follows that
\begin{eqnarray*}
\frac{ \|d\phi_{\tau_k(x)}(x)_{|E}\|}{e^{(\chi+ 2\ep) \tau_k(x)}} 
& \leq & \frac{\left(e^{a_3}\,e^{\psi^m(x)}\right)^{\taud s_0}}{e^{(\chi+ 2\ep) \tau_{ms_0}(x)}}
\leq  S_3\,\left(\frac{e^{\psi^m(x)}}{e^{(\tchi+ \ep) \frac{\tau_{m s_0}(x)}{\taud s_0}} } \right)^{\taud s_0}
\leq S_3 \, (D_\ep(x))^{\taud s_0} ,
\end{eqnarray*}
since $\frac{\tau_{ms_0}(x)}{s_0 \taud} = g^m(x)$. This is true for all $k \geq 0$, so
$C_{2\ep}(x) \leq S_3 \, (D_\ep(x))^{\taud s_0}$.
\endofproof

\bs

For any $p \in \R$ set
$$Q_{p}(\ep) = \{ x\in \ll \cap R  : D_\ep (x) \leq e^{p} \} ,$$
and, given $\hd \in (0,1]$,  denote by $\Xi_m(\ep,\hd)$ {\it the set of those $x\in \ll$ such that}
$$\sharp \, \left\{ k : 0 \leq k < m \: , \:  T^k(x) \notin Q_0(\ep) \right\} \geq \hd \, m .$$

The central result in this paper is the following.

\bs

\noindent
{\bf Theorem 4.4.} (Strong Exponential Tails)   {\it Let $\ep > 0$, $\delta \in (0, 1)$ and $\delta_0$ be
as in the beginning of Sect. 4, let the integer $s_0 = s_0 (\ep,\delta_0) \geq 1$ be chosen as above, and let
$T = \pp^{s_0}$. 
For every $\hd \in [\delta_0, 1]$ there 
exist constants  $C_2 = C_2(\ep, \hd) > 0$  and $c_2 = c_2(\ep,\hd) > 0$ such that 
\be
\mu\left( \Xi_m (\ep, \hd) \right) \leq C_2 e^{-c_2 m} .
\ee
In particular,
$$\mu\left( R \setminus \cup_{k=0}^{m-1} T^{-k} (Q_{0}(\ep)) \right) \leq C_2 e^{-c_2 m} $$
for every integer $m \geq 1$.}

\bs

We prove this theorem in the next section.
We will now derive Theorem 1.1 from it.

\bs

\noindent
{\it Proof of Theorem} 1.1. Recall that $\delta_0 = \frac{\delta}{2k_0}$ by (4.4). Take $\hd = \delta_0$.
Let $p_i \in \R$ be so that $S_i = e^{p_i}$ ($i = 1,2,3$),
where $S_1$, $S_2$ and $S_3$ are the constants from Lemmas 4.1 and 4.3. 

Given an integer $m \geq 0$, set $n = m s_0$ and  let $Y_n$ be {\it the set of those $x\in \ll $ such that}
$$\sharp \{ j : 0 \leq j < n \; , \; C_{2\ep}(\pp^j(x)) > e^{p_2 + p_3} \} \geq \delta_0  n .$$
Set $\Xi_m = \Xi_m(\ep, \delta_0)$ for brevity. We will now prove that 
\be
Y_n \subset \Xi_m .
\ee
Let $x\in Y_n$. Then $C_{2\ep}(\pp^j(x)) \leq e^{p_2+p_3}$ for less that $(1-\delta_0) n$ values
of $j = 0,1, \ldots, n-1$. By Lemma 4.3, 
$C_{2\ep}(\pp^j(x)) \leq S_3 (D_\ep(\pp^j(x)))^{\taud s_0} = e^{p_3} (D_\ep(\pp^j(x)))^{\taud s_0}$, so
$D_\ep(\pp^j(x)) \leq 1$ implies $C_{2\ep}(\pp^j(x)) \leq e^{p_3}$.
Assume for a moment that $x \notin \Xi_m$. Then $T^k(x) \in Q_0(\ep)$ for at least $(1-\delta_0)m$
values of $k = 0, 1 \ldots, m-1$. In other words $D_\ep(\pp^{ks_0}(x)) \leq 1$ for  at least $(1-\delta_0)m$
values of $k$, and so $C_{2\ep}(\pp^{ks_0}(x)) \leq e^{p_3}$ for at least $(1-\delta_0)m$
values of $k = 0, 1 \ldots, m-1$. For any such $k$, Lemma 4.1(b) implies that
$C_{2\ep}(\pp^{ks_0-j}(x)) \leq S_2 C_{2\ep}(\pp^{ks_0}(x)) \leq e^{p_2 + p_3}$ for all
$j = 0,1, \ldots, s_0-1$. Hence $C_{2\ep}(\pp^{r}(x)) \leq e^{p_2+p_3}$ for at least
$s_0 (1-\delta_0)m = (1-\delta_0)n$ values of $r = 0,1, \ldots,n-1$. In other words,
$C_{2\ep}(\pp^{r}(x)) > e^{p_2+p_3}$ for less than $\delta_0n$ values of $r = 0,1, \ldots,n-1$, which is 
a contradiction with $x\in Y_n$.

This proves (4.17). Combining the latter with (4.16) gives 
\be
\mu\left(Y_n \right) \leq C_2 e^{-c_2 n/s_0} .
\ee

Setting $\tY_n = \pi^{-1}(Y_n) = \{ y\in \tll : \pi(y) \in Y_n\}$, it follows from (4.1) that
\be
\m(\tY_n) = \frac{1}{\taud} \, \int_R \int_0^{\tau(y)} \chi_{\tY_n}(\phi_s(y))\, ds \; d\mu
= \frac{1}{\taud} \, \int_R \tau(y) \chi_{Y_n}(y)\; d\mu \leq \frac{\htau_0}{ \taud} \,  \mu(Y_n) .
\ee

Next, set  $k = [n/\tau_0]$, and 
let $X_k$ be {\it the set of those $y\in \tll$ such that}
$$\sharp \{ r : 0 \leq r < k \; , \; B_{2\ep}(f^r(y)) > e^{p_1 + p_2 + p_3} \} \geq \delta_0 k .$$
We will now prove that 
\be
X_k \subset \tY_n .
\ee
Given $y\in X_k$, set $x = \pi(y) \in \ll$; then $y = \phi_s(x)$ for some $s\in [0,\tau(x))$.
We have to show that $x\in Y_n$.
Now $y \in X_k$ means that $B_{2\ep}(f^r(x)) \leq e^{p_1 + p_2+p_3}$ for less then $(1-\delta_0) k$ values of 
$r = 0,1, \ldots k-1$. If $C_{2\ep}(\pp^j(x)) \leq e^{p_2+p_3}$ for some $j = 0,1, \ldots, n-1$, then for any
integer $r$ with 
\be
\tau_j(x) \leq r < \tau_j(x) + 1
\ee
we have $f^r(y) = \phi_t(\pp^j(x))$ for some
$t \in [-2,2]$, so by Lemma 4.1(a) we have $B_{2\ep}(f^r(y)) \leq S_1 C_{2\ep}(\pp^j(x)) \leq e^{p_1+p_2+p_3}$.
Assume for a moment that $x\notin Y_n$, i.e.  $C_{2\ep}(\pp^j(x)) \leq e^{p_2+p_3}$ for at least $(1-\delta_0) n$
values of $j = 0,1,\ldots,n-1$.  Notice that not more than $1/\tau_0$ values of $j$ will produce the same $r$
with (4.24), so there exist at least $(1-\delta_0) n/\tau_0 \geq (1-\delta_0)k$ values of $r$ with
$B_{2\ep}(f^r(y)) \leq e^{p_1+p_2+p_3}$. 
In other words, $B_{2\ep}(f^r(y)) > e^{p_1+p_2+p_3}$ for less than $\delta_0 k$
values of $r = 0,1, \ldots k-1$, a contradiction with $y \in X_k$.

This proves (4.20). Combining (4.20), (4.19) and (4.18) we get
$$\m(X_k) \leq \m(\tY_n) \leq \frac{\htau_0}{\taud} \, \mu(Y_n) \leq \frac{C_2 \htau_0}{\taud}\, e^{-c_2 n/s_0}
\leq \frac{C_2 \htau_0}{\taud}\, e^{-c_2 k\tau_0/s_0} = C_3 e^{-c_3 k} .$$

To complete the proof of Theorem 1.1, set $p_0 = p_1+p_2+p_3$.
Denote by $\Theta_k$ {\it the set of those $x\in \tll$ such that $\tR_{2\ep} (f^r(x)) > e^{p_0}$ for at least
$\delta k$ values of} $r = 0,1, \ldots, k-1$. We will prove an exponentially small estimate of $\m(\Theta_k)$.
Combining this with Proposition 3.1, in particular, using the left-hand-side inequality in (3.4), and (3.5) as well, will prove
Theorem 1.1.

Denote by  $Y^{+,i}_k$ the set of those $x\in \tll$ so that $\tR^+_{2\ep,i}(f^r(x)) > e^{p_0}$
for at least $\delta_0 k$ values of $r = 0,1, \ldots, k-1$. Define $Y^{-,i}_k$ in a similar way replacing $\tR^{+}_{2\ep,i}$ 
by $\tR^-_{2\ep,i}$. The above argument
shows that $\m(Y^{+,i}_k) \leq C_3 e^{-c_3 k}$ for all $i = 1, \ldots,k_0$.
Using Theorem 4.4 for the flow $\phi_{-t}$ and replacing the map $f = \phi_1$ by $f' = \phi_{-1}$, the 
above argument shows that  $\m(Y^{-,i}_k) \leq C_3 e^{-c_3 k}$. So, for the set 
$Y_k = \cup_{i=1}^{k_0} (Y^{+,i}_k \cup Y^{-,i}_k)$ we get $\m(Y_k) \leq C e^{-c k}$ for some constants
$C > 0$ and $c > 0$.

On the other hand, $\Theta_k \subset Y_k$. Indeed, let $x\in \Theta_k$. Then $\tR_{2\ep} (f^r(x)) > e^{p_0}$ for at least
$\delta k$ values of $r = 0,1, \ldots, k-1$. Since 
$\tR_{2\ep}(y) = \max_{1\leq i \leq k_0} \max\{\tR^+_{2\ep,i}(y), \tR^-_{2\ep,i}(y)\}$,
it follows that there exists $i = 1, \ldots, k_0$ such that either $\tR^+_{2\ep,i}(f^r(x)) > e^{p_0}$
for at least $\frac{\delta k}{2k_0} = \delta_0 k$ values of $r = 0,1, \ldots, k-1$ or 
$\tR^-_{2\ep,i}(f^r(x)) > e^{p_0}$ for at least $\frac{\delta k}{2k_0} = \delta_0 k$ values of $r = 0,1, \ldots, k-1$.
That is $x\in Y^{+,i}_k \cup Y^{-,i}_k \subset Y_k$. Hence $\Theta_k \subset Y_k$ and therefore
$\m(\Theta_k) \leq C e^{-c k}$.
\endofproof

\section{Proof of Theorem 4.4}
\setcounter{equation}{0}

Throughout we work under the assumptions in the beginning of Sect. 4 and will use the
notation from Sect. 4. Again, the constants $\ep > 0$, $\delta > 0$, $\delta_0 > 0$ and $s_0 = s_0(\ep, \delta_0)$ will be fixed as in 
Sect.4 and so will be the compact set $K_0$ and the functions $u$, $g$ and $\psi$.
As before $T = \pp^{s_0} : R \longrightarrow R$.

Setting
$$A_m(x) = \Aep_m(x) = \frac{e^{\psi^m(x)}}{e^{(\tchi + \ep) g^m(x)}} \quad , \quad x\in R ,$$
we have
$$D_\ep(x) = \sup_{m\geq 0} A_m(x) \geq 1 $$
for all $x\in \ll$. For such $x$ (modifying an idea in \cite{kn:Sim}) define
$$k_\ep(x) = \min \left\{ k \geq 0 : \psi^k(x) - (\tchi + \ep) g^k(x) = \log D_\ep(x) \right\} .$$
This is well-defined for $\mu$-almost all $x\in \ll$, since 
$\int_R g\, d\mu = 1$ and  $\int_R \psi\, d\mu = \tchi$, so
$$\lim_{k\to \infty} \frac{\psi^k(x)}{g^k(x)} = \lim_{k\to \infty} \frac{\psi^k(x)}{k} \cdot \frac{k}{g^k(x)} = \tchi $$
for $\mu$-almost all $x\in \ll$. Denote by $\ll'$ {\it the set of those $x\in \ll$ for which
$k_\ep(x)$ is well defined}; then $\mu(\ll') = 1$.

For any $n,k \geq 0$ we have $g^{n+k}(x) = g^k(x) + g^n(T^k(x))$ and also 
$\psi^{n+k}(x) = \psi^k(x) + \psi^n(T^k(x))$, therefore
\be
A_{n+k}(x) = A_k(x) \, A_n(T^k(x)) .
\ee
Using this inductively one derives that for every $x \in R$ and every integer $k \geq 1$ we have
\be
A_k(x) = A_1(x) A_1(T(x)) \ldots A_1 (T^{k-1}(x)) .
\ee

Next, consider the sets
$$\Omega_m(\ep) = \{ x\in \ll' : k_\ep(x) > m \} .$$

\ms

\noindent
{\bf Lemma 5.1.} {\it  There exist constants $C > 0$ and  $c  > 0$ such that 
\be
\mu \left(\Omega_m(\ep)\right) \leq C \, e^{- c m}
\ee
for all integers $m \geq 1$.}

\bs

\noindent
{\it Proof of Lemma} 5.1. Let $x \in \Omega_m(\ep)$. Then for $k = k_\ep(x)$ we have $k > m$ and
\be
\psi^k(x) - (\tchi + \ep) g^k(x) = \log D_\ep(x) \geq 0 .
\ee
Thus, $\psi^k(x)/k - (\tchi + \ep) g^k(x)/k \geq 0$, so
$$\left[\frac{\psi^k(x)}{k} - \tchi\right]  + \tchi\, \left[1 - \frac{g^k(x)}{k}\right]  
\geq \ep \, \frac{g^k(x)}{k} \geq \frac{\ep \tau_0}{\taud} .$$
Here we used $g^k(x) = \frac{\tau_{ks_0}(x)}{s_0 \taud} \geq \frac{k \tau_0}{\taud}$.

First assume $\chi = \chi_i \neq 0$; then (4.14) implies $|\tchi| \geq |\chi_i|/2 > 0$. 
Using the above, we either have
\be
\frac{\psi^k(x)}{k} - \tchi > \frac{\ep \tau_0}{2 \taud} ,
\ee
or 
\be
\left| 1 - \frac{g^k(x)}{k} \right| > \frac{\ep  \tau_0}{2 |\tchi| \taud} .
\ee
Let $\Delta'_k$ be the set of those $x\in \Omega_m(\ep)$ for which (5.5) holds, 
and let $\Delta''_k$ be the set of those $x\in \Omega_m(\ep)$ for which (5.6) holds.

We can now use the classical {\bf Large Deviation Principle} for the H\"older continuous functions $\psi$ and 
$g$ on $\hR$  (and the fact that $(\hR, T)$ is naturally isomorphic, up to a set of $\mu$-measure zero,  to a subshift of 
finite type\footnote{This natural isomorphism sends both $\psi$ and $g = \frac{1}{s_0\taud}\,\tau_{s_0}$
to H\"older continuous  functions with respect to the metric $d_\theta$ on the shift space.
Neither $\psi$ nor $g$ is cohomologous to a constant.} 
 -- see \cite{kn:OP}, \cite{kn:Ki} or \cite{kn:Y1}.  It follows from it that there exist constants\footnote{Notice that
the rate function involved depends on $T =\pp^{s_0}$ and therefore on $\ep$ and $\delta_0$. }
$C' > 0$ and $c' > 0$, independent of $k$ and $m$, such that 
$\mu(\Delta'_k) \leq C' e^{-c' k}$ and $\mu(\Delta''_k) \leq C' e^{-c' k}$. Since
$\Omega_m(\ep) \subset \cup_{k= m+1}^\infty (\Delta'_k \cup \Delta''_k)$, it follows that
$\mu(\Omega_m(\ep)) \leq C e^{-c m} $
for some constants $C, c > 0$ independent of $m$. 

When $\chi = \chi_i = 0$, (4.14) implies $0\leq \tchi < \frac{2C_1 \delta_0}{C_0}$, so (5.4) gives
$$\left[\frac{\psi^k(x)}{k} - \tchi\right]   \geq \ep \, \frac{g^k(x)}{k} + \tchi \, \frac{g^k(x)}{k} - \tchi 
\geq \ep \, \frac{g^k(x)}{k} -  \frac{8C_1 \delta_0}{C_0} \geq \frac{\ep \tau_0}{2\taud} ,$$
and then we can proceed as above.
\endofproof

\bs

The following lemma will be used later.

\bs

\noindent
{\bf Lemma 5.2.} {\it Let $x\in \ll'$ and let $k = k_\ep(x)$. Then:}

\ms

(a) {\it $T^{k}(x) \in Q_0(\ep)$.}

\ms

(b) {\it $k_\ep(x)$ coincides with the smallest integer $\ell \geq 0$ such that
$T^{\ell}(x) \in Q_0(\ep)$.}

\bs

\noindent
{\it Proof.} (a) From the definition of $k_\ep(x)$
it follows that $D_\ep(x) = A_{k_\ep(x)}(x) \geq A_{n+ k_\ep(x)}(x)$ for all $n \geq 0$. Using this and (5.1) with 
$k = k_\ep(x)$, we get $A_n(T^k(x)) \leq 1$ for all $n \geq 0$, so $D_\ep(T^k(x)) \leq 1$. 
Thus, $T^k(x) \in Q_0(\ep)$.

\ms

(b) Set $k = k_\ep(x)$. Assume that there exists $\ell = 0, 1, \ldots, k-1$ with $T^\ell(x) \in Q_0(\ep)$. Then
$A_n(T^\ell(x)) \leq 1$ for all $n \geq 0$, so by (5.1),
$A_{n+\ell}(x) = A_\ell(x) A_{n}(T^\ell(x)) \leq A_\ell(x)$
for all integers $n \geq 0$. In particular,
$A_k(x) = A_{(k-\ell) + \ell}(x) \leq A_{\ell}(x)$, which is a contradiction with the choice of $k = k_\ep(x)$.
\endofproof

\bs

The above already implies exponentially small tails.

\bs

\noindent
{\bf Lemma 5.3.} {\it There exist constants $C  > 0$ and $c  > 0$ such that
$$\mu\left(R\setminus   \cup_{k=0}^{m-1} T^{-k}(Q_0(\ep))\right) \leq C e^{-c m} $$
for all $m \geq 1$.}

\bs

\noindent
{\it Proof.} Let $x\in \ll'$ be such that $x\notin \Omega_m(\ep)$ for some $m \geq 1$, i.e. $k_\ep(x) < m$. 
Then by Lemma 5.2(a), $T^k(x) \in Q_0(\ep)$, so $x\in T^{-k} (Q_0(\ep))$ for some $k < m$. This shows that
$$\ll' \setminus \Omega_m(\ep) \subset \cup_{k=0}^{m-1} T^{-k}(Q_0(\ep)) \quad , \:\:\mbox{\rm i.e. } \quad
\ll' \setminus   \cup_{k=0}^{m-1} T^{-k}(Q_0(\ep)) \subset \Omega_m(\ep) ,$$
and by (5.3), 
$\mu\left(R\setminus   \cup_{k=0}^{m-1} T^{-k}(Q_0(\ep))\right) \leq C e^{-c m} ,$
which proves the Lemma.
\endofproof

\bs

Next, since the functions $\psi$ and $g$ are continuous on $\hR$ (with respect to the chosen metric $d_\theta$), the set
\be
F = \{ x\in \hR : \psi(x) - (\tchi + \ep/2) g(x) \leq 0 \}
\ee
is closed in $\hR$. Take $\eta = \frac{\delta_0}{C_0}$ and consider the open set
$$V_\eta = \{ x\in \hR : \psi(x) - (\tchi + \ep/2) g(x) < \eta \} .$$
Clearly, $F \subset V_\eta$, so there exists a H\"older continuous (with respect to $d_\theta$)
function $\varphi : \hR \longrightarrow [0,1]$
with $\varphi = 1$ on $F$ and $\varphi = 0$ on $\hR \setminus V_\eta$. {\bf Fix $\varphi$ with this property,}
and consider the H\"older continuous function 
$$h = (\psi - (\tchi + \ep/2) g) \, \varphi .$$
Notice that
\be
\left| \int_{\hR} h \, d\mu + \frac{\ep}{2} \right| < \frac{8 C_1 \delta_0}{C_0} .
\ee
Indeed, it follows from (4.10),  (4.11) and (4.14)  that for $x\in K_0$ we have
\begin{eqnarray*}
|(\psi(x) - (\tchi + \ep/2)g(x)) + \ep/2|
& \leq & |\psi(x) - \tchi| + |\tchi|\,  |1 - g(x)| + \frac{\ep}{2} |1 - g(x)|
\leq \frac{5 C_1 \delta_0}{C_0} .
\end{eqnarray*}
Since $\frac{5 C_1 \delta_0}{C_0} < \frac{\ep}{2}$ by (4.3), it follows that  $\psi(x) - (\tchi + \ep/2)g(x) < 0$ on $K_0$, 
so $K_0 \subset F$ and therefore $\varphi = 1$ on $K_0$. Combining the above with (4.9) gives
$$\left| \int_{\hR} h \, d\mu + \frac{\ep}{2} \right| \leq \left| \int_{K_0} (h +\ep/2)\, d\mu\right| 
+ (\|h\|_\infty + \ep/2) \mu (\hR\setminus K_0)
\leq \frac{5 C_1 \delta_0}{C_0} + \frac{3 C_1 \delta_0}{C_0} < \frac{8 C_1 \delta_0}{C_0} ,$$
since $\|h\|_\infty \leq \|\psi\|_\infty + C_1 \leq 3C_1$. This proves (5.8).

Apart from the above notice that 
\be
[\psi(y) - (\tchi + \ep/2) g(y)]\, \chi_F(y) \geq h(y) - \eta 
\ee
for all $y \in \ll'$. Here $\chi_F$ is the {\it characteristic function} of $F$. Indeed, if $y \in F$, then $\varphi(y) = 1$, so
$\psi(y) - (\tchi + \ep/2) g(y) = h(y)$ and (5.9) holds trivially. Let $x \notin F$; then
the left-hand-side of (5.9) is $0$. If moreover $x \notin V_\eta$, then $\varphi(y) = 0$,
so $h(y) = 0$, and (5.9) becomes $0 \geq -\eta$ which is obviously true. Finally, assume
$y \in V_\eta\setminus F$.  By the definitions of $F$ and $V_\eta$ we have 
$0 < \psi(y) - (\tchi + \ep/2) g(y) < \eta$. Multiplying this by $\varphi(y) \in  [0,1]$ we get
$0 \leq h(y) = (\psi(y) - (\tchi + \ep/2) g(y))\varphi (y) \leq \eta \varphi(y) \leq \eta$.
That is, $0 \geq h(y) - \eta$, so (5.9) holds again.


The following lemma proves Theorem 4.4. Recall the set $\Xi_m(\ep,\hd)$ defined just before Theorem 4.4.

\bs

\def\tvar{\tilde{\varphi}}

\noindent
{\bf Lemma 5.4.} {\it Let $\hd \in [\delta_0,1]$. There exist constants  
$C = C(\ep, \hd)\geq 1$ and $c = c(\ep, \hd) > 0$ such that
$\di \mu\left( \Xi_n(\ep,\hd) \right) \leq C e^{-c n} $
for any  integer $n \geq 1$.}

\bs

\noindent
{\it Proof.} Apart from the functions $A_m(x) = \Aep_m(x)$ defined in the beginning of Sect. 5 we will also use
the functions $\Aepp_m(x)$. Clearly
\be
\Aepp_m(x) = \frac{e^{\psi^m(x)}}{e^{(\tchi + \ep/2)g^m(x)}}
= \Aep_m(x) e^{\ep g^m(x)/2} \geq \Aep_m(x) e^{m \ep \tau_0/(2\taud)} \geq \Aep_m(x) 
\ee
for all $x\in \ll'$ and all integers $m \geq 0$. 


Notice that
$Q_0(\ep/2) \subset Q_0(\ep)$.
Indeed, if $y \in Q_0(\ep/2)$, then $\Aepp_m(y) \leq 1$ for all $m \geq 0$, and by (5.10),
$\Aep_m(y) \leq \Aepp_m(y) \leq 1$ for all $m \geq  0$, so $y \in Q_0(\ep)$.
Moreover, for any $y \in Q_0(\ep/2)$ we have $\Aepp_1(y) \leq 1$, that is $\psi(y) - (\tchi +\ep/2) g(y) \leq 0$, 
so $y\in F$, the set defined by (5.7). Thus, $Q_0(\ep/2) \subset F$. However in general $Q_0(\ep)$ is not
contained in $F$.

Next, setting $\delta_1 = \delta_0/C_0$, it follows from (5.3) that for any integer $n\geq 1$ we have
$\mu(\Omega_{[\delta_1 n]} (\ep) \leq C e^{-c (\delta_1 n-1)}$, so
$$\mu\left( \cup_{j=0}^n T^{-j} (\Omega_{[\delta_1 n]} (\ep)) \right) \leq (n+1) C e^{-c(\delta_1 n -1)} .$$
Setting
$$\Omn  =  \cup_{j=0}^n T^{-j} (\Omega_{[\delta_1 n]} (\ep)) ,$$
it follows from the above that there exist  constants $C_4 = C_4(\ep,\delta_0) > 0$ and $c_4 = c_4(\ep,\delta_0) > 0$
with
\be
\mu (\Omn) \leq C_4 e^{-c_4 n}
\ee
for all integers $n \geq 1$. Clearly, for every $x \in \ll' \setminus \Omn$ we have 
$k_{\ep}(T^j(x)) < \delta_1 n = \frac{\delta_0 n}{C_0}$ for all $j = 0,1, \ldots, n$.

\medskip

Next, {\bf fix for a moment an arbitrary}  $x\in \Xi_n \setminus \Omn$. We will now construct a sequence of points $x_i$ on the {\it orbit}
$$\oo_n(x) = \{ x, T(x), \ldots, T^{n-1}(x) \}$$
and certain integers $k_i$, $t_i$ using the set $Q_0(\ep)$. We will do the construction carefully and in all details,
although some of the details will not be used later. 

If $x\notin Q_0(\ep)$, set $x_0 = x$ and $t_0 = 0$. If $x\in Q_0(\ep)$ let $t_0 \geq 1$ be the largest 
integer such that $T^j(x) \in Q_0(\ep)$ for all $j =0,1,\ldots, t_0-1$. Then set $x_0 = T^{t_0}(x) \notin Q_0(\ep)$
and $k_0 = k_{\ep}(x_0) \geq 1$. Notice that $k_0 < \delta_1 n$.

By Lemma 5.2(a) we have $T^{k_0}(x_0) \in Q_0(\ep)$. 
Next, let $t_1 \geq 1$ be the largest  integer such that $T^{k_0+j}(x_0) \in Q_0(\ep)$ for all $j =0,1,\ldots, t_1-1$. 
Then set $x_1 = T^{k_0 + t_1}(x_0) \notin Q_0(\ep)$ and $k_1 = k_{\ep}(x_1) \geq 1$.

By induction we construct a sequence of points $x_0, x_1, \ldots, x_{s-1}$ and positive integers
$t_0, t_1, \ldots, t_{s-1}$ and $k_0, k_1, \ldots , k_{s-1}$ such that $x_{i+1} = T^{k_i+t_{i+1}}(x_i)$ for all
$i = 0,1, \ldots, s-2$, $k_i = k_{\ep}(x_i)$ for all $i = 0,1, \ldots, s-1$, and $t_{i+1} \geq 1$ is the maximal
integer such that $T^{k_i+j}(x_i) \in Q_0(\ep)$ for all $ j = 0,1 \ldots, t_{i+1}-1$. Thus, $x_i \notin Q_0(\ep)$,
so $k_i \geq 1$ for all $i = 0,1, \ldots, s-1$.

For $s$ we choose the maximal integer with
\be
\sum_{i=0}^{s-1} t_i + \sum_{i=0}^{s-1} k_i \leq n .
\ee
If we have equality in (5.12), set $t_s = 0$,  $k_s = 0$ and $x_s = T^{k_{s-1}}(x_{s-1})$.
If we have a strict inequality in (5.12), denote by $t_s$ the largest positive integer 
such that $T^{k_{s-1}+j}(x_{s-1}) \in Q_0(\ep)$ for all $ j = 0,1 \ldots, t_{s}-1$ and 
\be
\sum_{i=0}^{s} t_i + \sum_{i=0}^{s-1} k_i \leq n .
\ee
Set $x_s = T^{k_{s-1}+t_s}(x_{s-1})$. If there is an equality in (5.13), set $k_s = 0$.
If there is a strict inequality in (5.13), take $k_s \geq 0$ so that
\be
\sum_{i=0}^{s} t_i + \sum_{i=s}^{s} k_i = n .
\ee

It follows from the above construction that (5.14) always holds. Moreover, $k_i > 0$ for $i = 0,1,\ldots, s-1$,
while for $i = s$ we may have $k_s = 0$ and  $x_s \in Q_0(\ep)$. However, $k_s < \delta_1 n$ always holds, since
$x\notin \Omn$. Hence:
 
 (i) $T^{k_i+j}(x_i) \in Q_0(\ep)$ for all $j = 0,1 \ldots, t_{i+1}-1$, $i = 0,1, \ldots, s-1$,
 
 (ii) $k_i = k_{\ep}(x_i) \geq 1$ for all $ i = 0,1, \ldots, s-1$, and $0 \leq k_s \leq k_{\ep}(x_s) < \delta_1 n$.
 
 \ms
 
Thus, the above construction is such that the orbits $\oo_{k_i} (x_i)$ are disjoint, have no common points with $Q_0(\ep)$,  and
$$\oo_n(x) \setminus Q_0(\ep) = \cup_{i=0}^s \oo_{k_i} (x_i) .$$
 

 
 

In particular, it follows from it and the choice of $\psi$ in Lemma 4.2 that
$$|\log \Aepp_{k_s}(x_s)| = |\psi^{k_s}(x_s) - (\tchi + \ep/2) g^{k_s}(x_s) |
\leq  k_s(\|\psi\|_\infty + 2C_1) \leq 4 C_1 k_s <  4 C_1 \delta_1 n \leq   \frac{\ep \tau_0}{4 \taud} \hd n ,$$
using (4.3) and $\delta_1 = \delta_0/C_0 \leq \hd/C_0$.


Since $x\in \Xi_n = \Xi_n(\ep, \hd)$, it follows from (i) and (ii) that
\be
\sum_{i=0}^s k_i \geq  \hd n .
\ee

Next, it follows from (5.2) that for any $i = 0,1, \ldots,s$ we have
\be
\log \Aepp_{k_i}(x_i) = \sum_{q=0}^{k_i-1} \log \Aepp_1(T^q(x_i)) .
\ee
Similarly,
$$\psi^n(x) - (\tchi + \ep/2) g^n(x) = \log \Aepp_n(x) = \sum_{q=0}^{n-1} \log \Aepp_1(T^q(x)) ,$$
and it follows from the above construction that 
\begin{eqnarray}
\psi^n(x) - (\tchi + \ep/2) g^n(x) 
& = & \sum_{i=0}^s \log \Aepp_{k_i}(x_i) + \sum_{i=0}^{s-1} \log  \Aepp_{t_i}(T^{k_i}(x_i))\nonumber\\
& = & \sum_{i=0}^s \log \Aepp_{k_i}(x_i) + \sum_{q=0\atop {T^q(x) \in Q_0(\ep)}}^{n-1} \log  \Aepp_{1}(T^{q}(x)) .
\end{eqnarray}



By (ii) we have $\Aep_{k_i}(x_i) > 1$ for $i = 0, 1, \ldots, s-1$. This and (5.10) imply $\Aepp_{k_i}(x_i) > e^{\ep g^{k_i}(x_i)/2}$
and so $\log \Aepp_{k_i}(x_i) \geq \frac{\ep \tau_0}{2 \taud} k_i$ for all $i < s$. 
Combining this with the above, (5.15) and (5.17) give
$$\psi^n(x) - (\tchi + \ep/2) g^n(x) \geq \frac{\ep \tau_0}{2 \taud} \hd n - |\log \Aep_{k_s}(x_s) |
+ \sum_{q = 0\atop {T^q(x) \in Q_0(\ep)}}^{n-1} \Aepp_1(T^q(x))  \geq \frac{\ep \tau_0 \hd }{4 \taud}  n + (I) + (II) ,$$
where $(I)$ is the sum of the positive terms of the sum $\di \sum_{q = 0\atop {T^q(x) \in Q_0(\ep)}}^{n-1} \Aepp_1(T^q(x)) $
and $(II)$ is the sum of the other terms in this sum. Then $(I) \geq 0$, and, recalling the set $F$ defined by (5.7), 
$$(II) = \sum_{q = 0\atop {T^q(x) \in F\cap Q_0(\ep)}}^{n-1} \log \Aepp_1(T^q(x)) 
\geq \sum_{q = 0\atop {T^q(x) \in F}}^{n-1} \log \Aepp_1(T^q(x)).$$
(The terms appearing in the sum in the right-hand-side that do not appear in the sum in the middle are all non-positive.)
Hence
\begin{eqnarray*}
\psi^n(x) - (\tchi + \ep/2) g^n(x) 
& \geq & \frac{\ep \tau_0 \hd }{4 \taud}  n + \sum_{q = 0\atop {T^q(x) \in F}}^{n-1} \log \Aepp_1(T^q(x))\\
& =    & \frac{\ep \tau_0 \hd }{4 \taud}  n + \sum_{q = 0}^{n-1} [\psi (T^q(x)) - (\tchi + \ep/2) g(T^q(x))] \chi_F(T^q(x))\\
& =    & \frac{\ep \tau_0 \hd }{4 \taud}  n + \gamma^n(x) ,
\end{eqnarray*}
where 
$\gamma(y) = [\psi(y) - (\tchi + \ep/2)g(y)]\chi_F(y) .$
On the other hand (5.9) implies $\gamma^n(x) \geq h^n(x) - n \eta$. Using this and $\eta = \frac{\delta_0}{C_0}$,
it follows that
$$\psi^n(x) - (\tchi + \ep/2) g^n(x) \geq \frac{\ep \tau_0 \hd }{4 \taud}  n + h^n(x) - \frac{\delta_0}{C_0} n .$$
For the integral $\di H = \int_{\hR} h \, d\mu$ it follows from (5.8) that
$H + \frac{\ep}{2} > - \frac{8 C_1 \delta_0}{C_0}$.
So, we can now rewrite the above as
\begin{eqnarray*}
\left[\frac{\psi^n(x)}{n} - \tchi\right] + \tchi \, \left[1 - \frac{g^n(x)}{n}\right]
+ \frac{\ep}{2} \, \left[1 - \frac{g^n(x)}{n}\right] + \left[H - \frac{h^n(x)}{n}\right]
& \geq & H + \frac{\ep}{2} + \frac{\ep \tau_0 \hd }{4 \taud}  - \frac{\delta_0}{C_0} \\
& \geq & \frac{\ep \tau_0 \hd }{4 \taud}  - \frac{9 C_1 \delta_0}{C_0} \geq \frac{\ep \tau_0 \hd }{8 \taud}  ,
\end{eqnarray*}
since $\hd \geq \delta_0$ and $C_0 \geq \frac{100 C_1 \taud}{\ep \tau_0}$. Thus,
$$\left[\frac{\psi^n(x)}{n} - \tchi\right] + \tchi \, \left[1 - \frac{g^n(x)}{n}\right]
+ \frac{\ep}{2} \, \left[1 - \frac{g^n(x)}{n}\right] + \left[H - \frac{h^n(x)}{n}\right]
\geq \frac{\ep \tau_0 \hd }{8 \taud} .$$


As in the proof of Lemma 5.1, assume first that $\chi = \chi_i \neq 0$. Then
the above argument shows that $\Xi_n \setminus \Omn \subset \Xi'_n \cup \Xi''_n\cup \Xi_n'''\cup \Xi_n''''$, where
$\Xi'_n$ is the set of those $x\in \ll'\setminus \Omn$ such that
$\frac{\psi^n(x)}{n}  - \tchi > \frac{\ep \tau_0 \hd }{32 \taud}$, $\Xi''_n$
is the set of those $x\in \ll'\setminus \Omn$ such that
$\left| 1- \frac{g^n(x)}{n}\right| > \frac{\ep \tau_0 \hd }{32|\chi| \taud}$, $\Xi'''_n$
is the set of those $x\in \ll'\setminus \Omn$ such that
$\left| 1- \frac{g^n(x)}{n}\right| > \frac{\tau_0 \hd}{16 \taud}$, and $\Xi''''_n$
is the set of those $x\in \ll'\setminus \Omn$ such that
$\left| H- \frac{h^n(x)}{n}\right| > \frac{\ep \tau_0 \hd}{32 \taud}$.

As in the proof of Lemma 5.1, it now follows from the classical  Large Deviation Principle
(see \cite{kn:OP}, \cite{kn:Ki} or \cite{kn:Y1}) that there exist  constants $C' = C'(\ep, \hd) > 0$ and  $c' = c(\ep, \hd) > 0$
such that $\mu(\Xi_n \setminus \Omn) \leq C' e^{-c' n}$. 
By (5.11), $\mu(\Omn) \leq C_4 e^{-c_4 n}$, so  there exist constants $C \geq 1$
and $c > 0$ such that $\mu(\Xi_n) \leq C\, e^{-c n}$ for all $n$.

When $\chi = \chi_i = 0$, a slightly different argument (as in the proof of Lemma 5.1) shows that again 
there exist constants $C \geq 1$ and $c > 0$ such that $\mu(\Xi_n) \leq C\, e^{-c n}$ for all $n$.
\endofproof

\def\Abp{\overline{A}}
\def\tA{\widetilde{A}}

\section{Proof of Proposition 2.1}
\setcounter{equation}{0}

Throughout we assume that $\phi_t$ is a $C^2$ contact Anosov flow on $M$ with a $C^2$ invariant  contact
form $\omega$. Then the two-form $d\omega$ is $C^1$, so there exists a constant $C_0 > 0$ such that
\be
|d\omega_x(u,v)| \leq C_0 \|u\|\, \|v\| \quad , \quad u, v \in T_xM\:, \: x\in M ,
\ee
and for every $x\in M$ and every $u \in E^s(x)$ (or $u \in E^u(x)$)  with $\|u\| = 1$ there exists $v\in  E^u(x)$
(reps. $v \in E^s(x)$) with $\|v\| = 1$ such that $d\omega_x(u,v) \geq 1/C_0$. 

Fix for a moment $\ep > 0$. We will use the notation from Sect. 3.
Following general procedures (see Theorem S.2.10 in \cite{kn:KM}) for all $x\in \ll'_0$ and $i = 1, \ldots,k_0$
consider the inner product
$$\la u,v\ra'_{x,i} = \sum_{m\in \Z} \la df^m_x(u), df^m_x(v) \ra \, e^{-2 m \chi_i - 4\ep |m|} \quad , \quad 
u,v \in E_i(x) ,$$
where $\la \cdot, \cdot \ra$ is the inner product on $T_xM$ determined by the Riemann metric.
Then, given $u = u_1 + \ldots + u_{k_0}$, $v = v_1 + \ldots + v_{k_0} \in T_xM$ with $u_i,v_i \in E_i(x)$
for all $i$, define 
$$\la u, v\ra'_x = \sum_{i=1}^{k_0} \la u_i, v_i \ra'_{x,i} \quad ,\quad
\|u\|'_x = \sqrt{\la u, u\ra'_x} .$$
Notice that for $u\in E_i(x)$ we have 
\begin{eqnarray*}
\la u,u\ra'_{x,i} 
& =    & \sum_{m \geq 0} \|df^m_x(u)\|^2 \, e^{-2 m \chi_i - 4\ep m}
+ \sum_{m < 0} \|df^m_x(u)\|^2 \, e^{-2 m \chi_i + 4\ep m}\\
& \leq & \sum_{m \geq 0} (R^+_{\ep}(x) e^{m\chi_i+ m\ep}\|u\|)^2 \, e^{-2 m \chi_i - 4\ep m}
+ \sum_{m < 0} (R^-_{\ep}(x) e^{m\chi_i- m\ep}\|u\|) ^2 \, e^{-2 m \chi_i + 4\ep m}\\
& =    & (R_{\ep}(x) \|u\|)^2 \left[\sum_{m\geq 0} e^{-2 m\ep} + \sum_{m< 0} e^{2 m\ep} \right]
\leq \frac{1+ e^{-2\ep}}{1- e^{-2\ep}} (R_{\ep}(x) \|u\|)^2 \leq \frac{2}{\ep} (R_{\ep}(x) \|u\|)^2 .
\end{eqnarray*}
Thus, $\|u\|'_{x,i} \leq \sqrt{2/\ep} \; R_{\ep}(x) \|u\|$.

Since the subspaces $E_i(x)$ are mutually orthogonal with respect to $\la \cdot, \cdot \ra'_x$, 
it follows that, setting  $K_0 = K_0(\ep) = 2 \sqrt{k_0/\ep}$, we have
\be
\|u\|'_x \leq K_0(\ep) \, R_{\ep}(x) \|u\| \quad, \quad u \in T_xM \:, \: x\in \ll'_0 .
\ee
In particular, if $u = u_1 +\ldots + u_{k_0}$ with $u_i \in E_i(x)$ for all $i$, then 
\be
\|u_i\| \leq \|u_i\|'_{x,i} \leq \|u\|'_x \leq K_0(\ep) \, R_{\ep}(x) \, \|u\|  \quad , \quad
i = 1, \ldots, k_0 .
\ee

Next, for every $i = 1, \ldots,k_0$, $x\in \ll'_0$ and  $\ep > 0$ set
$$\Abp_i(x,\ep) = \sup \left\{ \frac{\| df_x^{n+k}(v)\| \; e^{-(\chi_i+\ep)n - \ep |k|}}{\|df_x^k (v)\|} \; : \: n\geq 0 , k \in \Z,
v\in E_i\setminus \{0\} \right\} ,$$
$$\tA_i(x, \ep) = \sup \left\{ \frac{\| df_x^{n+k}(v)\| \; e^{-(\chi_i-\ep)n - \ep |k|}}{\|df_x^k (v)\|} \; : \: n\leq 0 , k \in \Z,
v\in E_i\setminus \{0\} \right\} ,$$
and
$$A_\ep(x) = \max_{1\leq i \leq k_0 } \: \max \{ \Abp_i(x,\ep) , \tA_i(x, \ep) \} .$$
It follows from the argument in the proof of Proposition 1 in \cite{kn:FHY}
that $\Abp_i(x,\ep)$ and
$\tA_i(x, \ep)$ are finite for all $i$, so $A_\ep(x)$ is well-defined for all $x\in \ll'_0$, and moreover
\be
\frac{\|df_x^n(v)\|}{e^{n\chi_i}} 
\leq A _\ep(x)\, e^{|n|\ep} \|v\| \quad , \quad v\in E_i(x)  \;, \; n \in \Z ,
\ee
for all $x\in \ll'_0$ and all $i =1, \ldots,k_0$. It now follows from (3.1) and (3.2) that
$$R_\ep(x) \leq A_\ep(x)$$
for all $\ep > 0$ and all $x\in \ll'_0$.  It is easy to see that $\Abp_i(x,\ep)$ and
$\tA_i(x, \ep)$ are $\ep$-slowly varying functions for all $i$, so $A_\ep(x)$ is also an
$\ep$-slowly varying function.

Next, set $K_1(\ep) = C_0^2 K_0(\ep)$. We will prove that
\be
\frac{1}{K_1(\ep) \, (R_\ep(x))^2 e^{m\ep}} \leq \frac{\|df_x^m (v)\|}{e^{m \chi_i}} \quad, \quad m \geq 0\;,\;
v \in E_i(x)\;,\; \|v\| = 1,
\ee
for all $x\in \ll'_0$ and $i = 1, \ldots,k_0$. Given $v\in E_i(x)$ with $\|v\| = 1$, there exists
$w\in T_xM$ with $\|w\| = 1$ such that $d\omega_x(v,w) \geq 1/C_0$. Let $w = w_1 + \ldots + w_{k_0}$,
where $w_k \in E_k(x)$ for all $k$. Since the flow is contact, there exists 
$j$ with $\chi_j = -\chi_i$. It is now easy to see that $d\omega_x(v,w) = d\omega_x(v, w_j)$.
Indeed, setting $x_m = f^m(x)$, for $\chi_k < \chi_j = -\chi_i$ we have $\chi_k < \chi_j - 3\ep$ by the choice of $\ep$, 
so using (6.1) we get 
\begin{eqnarray*}
|d\omega_x(v,w_k)| 
& =     &  |d\omega_{x_m}(df^m_x (v), df^m_x (w_k))| \leq  C_0 \|df^m_x (v)\| \,\|df^m_x (w_k)\|\\
& \leq & C_0 R_{\ep}^2(x) \|v\|\, \| w_k\|\, e^{m \chi_i + m\ep} e^{m \chi_k + m \ep}  \to 0
\end{eqnarray*} 
as $m \to \infty$. Thus, $d\omega_x(v,w_k) = 0$. In a similar way one deals with the case $\chi_k > \chi_j$,
considering $m \to -\infty$.

Hence for any $m \geq 0$ we have
$$\frac{1}{C_0} \leq d\omega_x(v,w_j) = d\omega_{f^m(x)} (df^m_x(v), df^m_x(w_j))
\leq C_0 \| df^m_x(v)\| \, \|df^m_x(w_j)\| .$$
This, (3.1) and (6.3) imply 
\begin{eqnarray*}
\frac{\| df^m_x(v)\|}{e^{m \chi_i}} 
& \geq & \frac{1}{C_0^2 \|df^m_x(w_j)\| e^{m \chi_i}} \geq \frac{1}{C_0^2 R^+_{\ep}(x) e^{m\chi_j + m \ep} \|w_j\| e^{m \chi_i}}\\
& \geq & \frac{1}{C_0^2 R^+_{\ep}(x) e^{m \ep} K_0(\ep) R_\ep(x) \|w\|}
\geq  \frac{1}{C_0^2 K_0(\ep) (R_{\ep}(x))^2 e^{m \ep} } ,
\end{eqnarray*}
which proves (6.5). In a similar way one shows that
\be
\frac{1}{K_1(\ep) \, (R_\ep(x))^2 e^{m\ep}} \leq \frac{\|df_x^{-m} (v)\|}{e^{-m \chi_i}} \quad, \quad m \geq 0\;,\;
v \in E_i(x)\;,\; \|v\| = 1,
\ee
for all $x\in \ll'_0$ and $i = 1, \ldots,k_0$. 

Next, we will prove that
\be
\Abp_i(x,\ep) \leq K_1(\ep/2) (R_{\ep/2}(x))^3
\ee
for all $x\in \ll'_0$ and $i = 1, \ldots, k_0$. Fix $x$ and $i$ for a moment. Given $v\in E_i(x) \setminus \{0\}$,
$n \geq 0$ and $k \in \Z$,  we  will consider 3 cases.

\ms

{\bf Case 1.} $k \geq 0$; then $n+k \geq 0$, too. Using (3.1) and (6.5) with $\ep$ replaced by $\ep/2$, we get
\begin{eqnarray*}
\frac{\| df_x^{n+k}(v)\| \; e^{-(\chi_i+\ep)n - \ep |k|}}{\|df_x^k (v)\|} 
& \leq & \frac{R_{\ep/2}(x) e^{ (n+k) (\chi_i + \ep/2)}\|v\| 
\; e^{-(\chi_i+\ep)n - \ep |k|}}{\frac{e^{k \chi_i} \|v\|}{K_1(\ep/2) \, (R_{\ep/2}(x))^2 e^{k\ep/2}} } \\
& =    & K_1(\ep/2) (R_{\ep/2}(x))^3 e^{ (n+k) \ep/2} e^{- n\ep - k \ep } e^{k\ep/2}\\
& \leq & K_1(\ep/2) (R_{\ep/2}(x))^3 .
\end{eqnarray*}


{\bf Case 2.} $k < 0$ and $n+k \geq 0$. Then (3.1) and (6.6) imply
\begin{eqnarray*}
\frac{\| df_x^{n+k}(v)\| \; e^{-(\chi_i+\ep)n - \ep |k|}}{\|df_x^k (v)\|} 
& \leq & \frac{R_{\ep/2}(x) e^{ (n+k) (\chi_i + \ep/2)}\|v\| 
\; e^{-(\chi_i+\ep)n + \ep k}}{\frac{e^{k \chi_i} \|v\|}{K_1(\ep/2) \, (R_{\ep/2}(x))^2 e^{-k\ep/2}} } \\
& =    & K_1(\ep/2) (R_{\ep/2}(x))^3 e^{ (n+k) \ep/2} e^{- n\ep + k \ep } e^{-k\ep/2}\\
& \leq & K_1(\ep/2) (R_{\ep/2}(x))^3 .
\end{eqnarray*}

{\bf Case 3.} $k < 0$ and $n+k < 0$. Then (3.2) and (6.6) imply
\begin{eqnarray*}
\frac{\| df_x^{n+k}(v)\| \; e^{-(\chi_i+\ep)n - \ep |k|}}{\|df_x^k (v)\|} 
& \leq & \frac{R_{\ep/2}(x) e^{ (n+k) (\chi_i - \ep/2)}\|v\| 
\; e^{-(\chi_i+\ep)n + \ep k}}{\frac{e^{k \chi_i} \|v\|}{K_1(\ep/2) \, (R_{\ep/2}(x))^2 e^{-k\ep/2}} } \\
& =    & K_1(\ep/2) (R_{\ep/2}(x))^3 e^{ -(n+k) \ep/2} e^{- n\ep + k \ep } e^{-k\ep/2}\\
& \leq & K_1(\ep/2) (R_{\ep/2}(x))^3 .
\end{eqnarray*}

This proves (6.7). In a similar way we will now prove
\be
\tA_i(x,\ep) \leq K_1(\ep/2) (R_{\ep/2}(x))^3
\ee
for all $x\in \ll'_0$ and $i = 1, \ldots, k_0$. Fix $x$ and $i$ for a moment and let $v\in E_i(x) \setminus \{0\}$,
$n \leq 0$ and $k \in \Z$. Again  we  will consider 3 cases.

\ms

{\bf Case 4.} $k \geq 0$ and $n+k \geq 0$. Using (3.1) and (6.5), we get
\begin{eqnarray*}
\frac{\| df_x^{n+k}(v)\| \; e^{-(\chi_i - \ep)n - \ep |k|}}{\|df_x^k (v)\|} 
& \leq & \frac{R_{\ep/2}(x) e^{ (n+k) (\chi_i + \ep/2)}\|v\| 
\; e^{-(\chi_i-\ep)n - \ep k}}{\frac{e^{k \chi_i} \|v\|}{K_1(\ep/2) \, (R_{\ep/2}(x))^2 e^{k\ep/2}} } \\
& =    & K_1(\ep/2) (R_{\ep/2}(x))^3 e^{ (n+k) \ep/2} e^{n\ep - k \ep } e^{k\ep/2}\\
& \leq & K_1(\ep/2) (R_{\ep/2}(x))^3 .
\end{eqnarray*}


{\bf Case 5.} $k \geq 0$ and $n+k < 0$. Then 
\begin{eqnarray*}
\frac{\| df_x^{n+k}(v)\| \; e^{-(\chi_i-\ep)n - \ep |k|}}{\|df_x^k (v)\|} 
& \leq & \frac{R_{\ep/2}(x) e^{ (n+k) (\chi_i - \ep/2)}\|v\| 
\; e^{-(\chi_i-\ep)n - \ep k}}{\frac{e^{k \chi_i} \|v\|}{K_1(\ep/2) \, (R_{\ep/2}(x))^2 e^{k\ep/2}} } \\
& =    & K_1(\ep/2) (R_{\ep/2}(x))^3 e^{ -(n+k) \ep/2} e^{n\ep - k \ep } e^{k\ep/2}\\
& \leq & K_1(\ep/2) (R_{\ep/2}(x))^3 .
\end{eqnarray*}

{\bf Case 6.} $k < 0$, and so $n+k < 0$. Then 
\begin{eqnarray*}
\frac{\| df_x^{n+k}(v)\| \; e^{-(\chi_i-\ep)n - \ep |k|}}{\|df_x^k (v)\|} 
& \leq & \frac{R_{\ep/2}(x) e^{ (n+k) (\chi_i - \ep/2)}\|v\| 
\; e^{-(\chi_i-\ep)n + \ep k}}{\frac{e^{k \chi_i} \|v\|}{K_1(\ep/2) \, (R_{\ep/2}(x))^2 e^{-k\ep/2}} } \\
& =    & K_1(\ep/2) (R_{\ep/2}(x))^3 e^{ -(n+k) \ep/2} e^{ n\ep + k \ep } e^{-k\ep/2}\\
& \leq & K_1(\ep/2) (R_{\ep/2}(x))^3 .
\end{eqnarray*}

This proves (6.8). Combining the latter with (6.7) gives $A_\ep(x) \leq K_1(\ep/2) (R_{\ep/2}(x))^3$, and therefore
$$K_2(\ep) (A_{2\ep}(x))^{1/3} \leq R_\ep(x) \leq A_\ep(x) \quad ,\quad x\in \ll'_0 ,$$
where $K_2(\ep) = \frac{\ep^{1/6}}{C_0^{2/3} (2k_0)^{1/6}}$.
\endofproof

\bs

\noindent
{\it Proof of } (3.5). Given $i = 1, \ldots,k_0$ and $x \in \ll'_0$, there exist an integer $n \geq 0$ and $v \in \Ei$
with $\|v\| = 1$ such that
$\tR^+_{\ep,i}(x) = \frac{\|df^n_x (v)\|}{e^{(\chi_i+\ep)}}$.
Let $v = \sum_{k=1}^i v_k$ with $v_k \in E_k(x)$ for all $k$. Then by (6.3),
$\|v_k\| \leq K_0(\ep) R_\ep(x) \|v\| \leq K_0(\ep) R_\ep(x)$ for all $k$. Using this and $\chi_i \geq \chi_k$ for 
$k\leq i$ yields
\begin{eqnarray*}
\tR^+_{\ep,i}(x) 
& =    &  \frac{\|df^n_x (v)\|}{e^{(\chi_i+\ep)n}} \leq  \frac{\sum_{k=1}^i \|df^n_x (v_k)\|}{e^{(\chi_i+\ep)n}}
\leq \sum_{k=1}^i  \frac{\|df^n_x (v_k)\|}{e^{(\chi_k+\ep) n}}\\
& \leq & \sum_{k=1}^i  \frac{\|(df^n_x)_{| E_k}\| \, \|v_k\|}{e^{(\chi_k+\ep)n}}
\leq \sum_{k=1}^i R^+_{\ep,k}(x)\;  K_0(\ep) R_\ep(x) \leq k_0 K_0(\ep)\, (R_\ep(x))^2 .
\end{eqnarray*}
This is true for all $i$, so $\tR^+_{\ep}(x) \leq k_0 K_0(\ep)\, (R_\ep(x))^2$. In a similar way, considering 
the flow $\phi_{-t}$ and the map $f^{-1} = \phi_{-1}$ one shows that $\tR^-_{\ep}(x) \leq k_0 K_0(\ep)\, (R_\ep(x))^2$.
This proves (3.5).
\endofproof

\bs

{\sc University of Western Australia, Crawley WA 6009, Australia}

{\sc\it E-mail address:} luchezar.stoyanov@uwa.edu.au

\end{document}